\documentclass[11pt]{article}
\usepackage{amsmath} \usepackage{amssymb} \usepackage{latexsym} \usepackage{enumitem} \usepackage{mathrsfs} \usepackage{comment}
\usepackage{color}

\usepackage[colorlinks=true,urlcolor=blue,
citecolor=red,linkcolor=blue,linktocpage,pdfpagelabels,bookmarksnumbered,bookmarksopen]{hyperref}

\newtheorem{assumption}{Assumption}

\def\qed{ \ \vrule width.2cm height.2cm depth0cm\smallskip}
\newenvironment{proof}{\noindent {\bf Proof.\/}}{$\qed$\vskip 0.1in}

\newcommand{\ol}{\overline}
\newcommand{\ul}{\underline}

\newcommand{\ba}{\begin{array}}
\newcommand{\ea}{\end{array}}
\newcommand{\be}{\begin{equation}}
\newcommand{\ee}{\end{equation}}
\newcommand{\bea}{\begin{eqnarray}}
\newcommand{\eea}{\end{eqnarray}}
\newcommand{\beaa}{\begin{eqnarray*}}
\newcommand{\eeaa}{\end{eqnarray*}}

\def\dbC{\mathbb{C}}

\def\dbE{\mathbb{E}}
\def\dbF{\mathbb{F}}

\def\dbL{\mathbb{L}}

\def\dbP{\mathbb{P}}
\def\dbR{\mathbb{R}}

\def\a{\alpha}
\def\b{\beta}
\def\g{\gamma}
\def\d{\delta}
\def\e{\varepsilon}

\def\k{\kappa}
\def\l{\lambda}

\def\si{\sigma}

\def\f{\varphi}
\def\th{\theta}
\def\o{\omega}
\def\h{\widehat}
\def\G{\Gamma}

\def\Th{\Theta}

\def\O{\Omega}
\def\cC{{\cal C}}

\def\cF{{\cal F}}
\def\cG{{\cal G}}

\def\cL{{\cal L}}

\def\cN{{\cal N}}

\def\cP{{\cal P}}

\def\no{\noindent}

\def\ss{\smallskip}
\def\ms{\medskip}

\def\q{\quad}
\def\qq{\qquad}

\def\pa{\partial}
\def\cd{\cdot}
\def\cds{\cdots}

\def\td{\nabla}

\def\tr{\hbox{\rm tr}}

\def\qed{ \hfill \vrule width.25cm height.25cm depth0cm\smallskip}

\newcommand{\basa}{\begin{assumption}}
\newcommand{\easa}{\end{assumption}}

\newcommand{\bas}{\begin{assum}}
\newcommand{\eas}{\end{assum}}

\def\pa{\partial}
\def\h{\widehat}

 \def\cd{\cdot}
\def\cds{\cdots}

\def\tr{\hbox{\rm tr$\,$}}

\def\dis{\displaystyle}

\def\wh{\widehat}

\def\1{{\bf 1}}

\def\:{\!:\!}
\def\reff{\eqref}
\def \proof{{\noindent \bf Proof.\quad}}

at 9pt

\definecolor{alp}{rgb}{0.0, 0.5, 0.0}

\def\R{\mathbb{R}}

\newtheorem{thm}{Theorem}[section]
\newtheorem{lem}[thm]{Lemma}

\newtheorem{prop}[thm]{Proposition}
\newtheorem{rem}[thm]{Remark}
\newtheorem{eg}[thm]{Example}
\newtheorem{defn}[thm]{Definition}
\newtheorem{assum}[thm]{Assumption}

\begin{document}

\title{\bf Mean Field Game Master Equations with Anti-monotonicity Conditions } 
\author{Chenchen Mou\thanks{\noindent  Dept. of Math.,
City University of Hong Kong. E-mail: \href{mailto:chencmou@cityu.edu.hk}{chencmou@cityu.edu.hk}. This author is supported in part by CityU Start-up Grant 7200684 and Hong Kong RGC Grant ECS 9048215.} ~ and ~ Jianfeng Zhang\thanks{\noindent  Dept. of Math., 
University of Southern California. E-mail:
\href{mailto:jianfenz@usc.edu}{jianfenz@usc.edu}. This author is supported in part by NSF grants DMS-1908665 and DMS-2205972.
}  
}
\date{}
\maketitle

\begin{abstract} It is well known that the monotonicity condition, either in Lasry-Lions sense or in displacement sense, is crucial for the global well-posedness of  mean field game master equations, as well as for the uniqueness of mean field equilibria and solutions to mean field game systems. In the literature, the monotonicity conditions are always taken in a fixed direction. In this paper we propose a new type of monotonicity condition in the opposite direction, which we call the anti-monotonicity condition, and establish the global well-posedness for mean field game master equations with nonseparable Hamiltonians. Our anti-monotonicity condition allows our data to violate both the Lasry-Lions monotonicity and the displacement monotonicity conditions.
\end{abstract}

\no{\bf Keywords.}  Master equation, mean field games, Lasry-Lions monotonicity, displacement monotonicity, anti-monotonicity.

\ms
\no{\it 2020 AMS Mathematics subject classification:}  35R15, 49N80, 49Q22, 60H30, 91A16, 93E20

\vfill\eject


\section{Introduction}
\label{sect-Introduction}
\setcounter{equation}{0} 

In this paper we consider the following second order master equation, arising from mean field games with common noise, with terminal condition $V(T,x,\mu) = G(x,\mu)$:
\bea
\label{master}
\left.\ba{lll}
\dis\cL V(t,x,\mu):= -\pa_t V -\frac{\h\b^2}{2} \tr(\pa_{xx} V) + H(x,\mu,\partial_x V)  - \cN V =0,  \q\mbox{where}\\
\dis \cN V(t,x,\mu):= \tr\bigg(\bar{\tilde  \dbE}\Big[\frac{\h\b^2}{2} \pa_{\tilde x} \pa_\mu V(t,x, \mu, \tilde \xi) - \pa_\mu V(t, x, \mu, \tilde \xi)(\pa_pH)^\top(\tilde \xi,\mu, \pa_x V(t, \tilde \xi, \mu))
 \\
 \dis\qq+\b^2 \pa_x\pa_\mu V(t,x,\mu,\tilde \xi)+\frac{\b^2}{2}\pa_{\mu\mu}V(t,x,\mu,\bar\xi,\tilde\xi)\Big]\bigg),\q   (t,x, \mu)\in [0,T)\times\R^d\times\cP_2(\R^d).
 \ea\right.
\eea
Here $\b\geq 0$ is a constant, $\h\b^2:= 1+\b^2$,  $\partial_t,\partial_x,\partial_{xx}$ are standard temporal and spatial derivatives, $\partial_{\mu},\partial_{\mu\mu}$ are $W_2$-Wasserstein derivatives, $\tilde \xi$ and $\bar\xi$ are independent random variables with the same law $\mu$ and $\bar {\tilde \dbE}$ is the expectation with respect to their joint law. The theory of Mean Field Games (MFGs, for short), initiated independently by Caines-Huang-Malham\'e \cite{HCM06} and Lasry-Lions \cite{LL07a}, studies the asymptotic behavior of stochastic differential games with a large number of players interacting in certain symmetric way. We refer to Lions \cite{Lions}, Cardaliaguet \cite{Cardaliaguet}, Bensoussan-Frehse-Yam \cite{BFY}, Carmona-Delarue \cite{CD1,CD2} and Cardaliaguet-Porretta \cite{CP} for a comprehensive exposition of the subject. First introduced by Lions \cite{Lions}, the master equation characterizes the value of the MFG, provided there is a unique mean field equilibrium. Roughly speaking, it plays the role of the HJB equation in the stochastic control theory.

The master equation \reff{master} admits a unique local (in time) classical solution when the data $H$ and $G$ are sufficiently smooth, see e.g. Gangbo-Swiech \cite{GangboS2015}, Bensoussan-Yam \cite{BY}, Mayorga \cite{Mayorga}, Carmona-Delarue \cite{CD2} and Cardaliaguet-Cirant-Porretta \cite{CCP}. In particular, \cite{CCP} studied the local well-posedness of the master equations not only for MFGs involving homogeneous minor players but also for MFGs with a major player. It is much more challenging to obtain a global classical solution, we refer to Buckdahn-Li-Peng-Rainer \cite{BLPR}, Chassagneux-Crisan-Delarue \cite{CCD}, Cardaliaguet-Delarue-Lasry-Lions \cite{CDLL}, Carmona-Delarue \cite{CD2}, Gangbo-Meszaros-Mou-Zhang \cite{GMMZ} and, in the realm of potential MFGs,  Bensoussan-Graber-Yam \cite{BGY1,BGY2}, Gangbo-Meszaros \cite{GM}. We also refer to Mou-Zhang \cite{MZ2}, Bertucci \cite{B1}, and Cardaliaguet-Souganidis \cite{CarSou} for global weak solutions which require much weaker regularity on the data, and Bayraktar-Cohen \cite{BC}, Bertucci-Lasry-Lions \cite{BLL}, Cecchin-Delarue \cite{CecchinDelarue}, Bertucci \cite{Bertucci0} for classical or weak solutions of finite state mean field game master equations. All the above global well-posedness results, with the exception \cite{BLPR} that considers linear master equations and thus no control or game is involved, require certain monotonicity condition, which we explain next.

One typical condition, extensively used in the literature \cite{BC,Bertucci0,B1,BLL,CDLL,CarSou,CD2,CCD,MZ2}, is the well-known Lasry-Lions monotonicity condition: for a function $G:\mathbb R^d\times\mathcal{P}_2(\mathbb R^d)\to\mathbb R$, 
\bea
\label{LLmon}
\mathbb E\Big[G(\xi_1,\mathcal{L}_{\xi_1})+G(\xi_2,\mathcal{L}_{\xi_2})-G(\xi_1,\cL_{\xi_2})-G(\xi_2,\mathcal{L}_{\xi_1})\Big]\ge 0,
\eea
for any square integrable random variables $\xi_1, \xi_2$.  Another type of monotonicity condition,originating in  Ahuja \cite{Ahuja} and was later sparsely used in the literature,  see Ahuja-Ren-Yang \cite{ARY} and \cite{BGY1, BGY2, GM, GMMZ}, is the displacement (or weak) monotonicity, 
\bea
\label{dismon}
\mathbb E\Big[[\pa_xG(\xi_1,\mathcal{L}_{\xi_1})-\pa_xG(\xi_2,\mathcal{L}_{\xi_2})][\xi_1-\xi_2]\Big]\ge 0.
\eea
When $G$ is regular enough with bounded $\pa_{xx}G,\pa_{x\mu}G$, \reff{LLmon} and \reff{dismon} are equivalent to the following inequalities, respectively: for all square integrable random variables $\xi,\eta$,
\bea
\label{monsmooth}
 \tilde \dbE \big[\big \langle \partial_{x \mu} G(\xi, \cL_\xi, \tilde \xi) \tilde \eta, \eta\big \rangle\big]\ge 0,\q \tilde \dbE \big[\big \langle \partial_{x \mu} G(\xi, \cL_\xi, \tilde \xi) \tilde \eta, \eta\big \rangle\big]+\mathbb E\big[\langle \partial_{xx} G(\xi, \cL_\xi) \eta, \eta\big \rangle\big]\ge 0,
\eea
where $(\tilde \xi, \tilde \eta)$ is an independent copy of $(\xi, \eta)$. The monotonicity conditions are crucial for the uniqueness of the Nash equilibria  of MFGs and thus the well-posedness of their master equations.  

When none of the monotonicity conditions holds, the MFG could have multiple equilibria, see e.g. Foguen Tchuendom \cite{F}, Cecchin-Dai Pra-Fisher-Pelino \cite{CDFP}, Bayraktar-Zhang \cite{BZ}. In this case, one approach is to consider a special type of equilibria,
 see e.g.  \cite{CDFP}, Delarue-Foguen Tchuendom \cite{DF}, Cecchin-Delarue \cite{CecchinDelarue}, Bayraktar-Cecchin-Cohen-Delarue \cite{BCCD1,BCCD2}. A larger literature is on the possible convergence of the equilibria for the $N$-player game, which is quite often unique because the corresponding Nash system is non-degenerate due to the presence of the individual noises, to the mean field equilibria (which may or may not be unique), see, e.g., \cite{CDLL, CD2, MZ2}, Delarue-Lacker-Ramanan \cite{DLR1,DLR2}, Djete \cite{Djete}, Lacker \cite{Lacker1,Lacker2,Lacker3,Lacker4}, Lacker-Flem \cite{LackerFlem}, Nuts-San Martin-Tan \cite{NSMT}. Finally, we note that Iseri-Zhang \cite{IZ} takes a quite different approach by investigating the set of game values over all mean field equilibria  and establishes the dynamic programming principle and the convergence from the $N$-player game to the MFG. 

We emphasize that the two inequalities in \reff{monsmooth} share the same direction. Our goal of this paper is to propose a new type of monotonicity condition in the opposite direction, which we call anti-monotonicity condition, and establish the global well-posedness 
for the master equation \reff{master}, with possibly nonseparable Hamiltonian $H$. We remark that the mean field equilibrium is a fixed point, and the monotonicity conditions \reff{monsmooth} were used to ensure the uniqueness of the fixed point. To motivate our anti-monotonicity condition, let us use a  very simple example to illustrate the idea. Suppose that $f:\mathbb R^1\to\mathbb R^1$ is a continuously differentiable function and we are interested in its fixed point $x^*$: $f(x^*) = x^*$. When $f$ is decreasing, i.e., $f'\leq 0$,  clearly $f$ admits a unique fixed point $x^*$. When $f$ is increasing, in general neither the existence nor the uniqueness of $x^*$ is guaranteed. However, if $f$ is sufficiently monotone, in the sense that $f'\ge 1+\e$ for some $\epsilon>0$, then again $f$ has a unique fixed point $x^*$. While in complete different contexts, our conditions follow the same spirit. Roughly speaking, the standard monotonicity conditions \reff{monsmooth} correspond to the case that $f$ is decreasing, while our new anti-monotonicity condition corresponds to the case $f$ is increasing, and for the same reason we will need to require our data to be sufficiently anti-monotone in appropriate sense.

To be precise, our anti-monotonicity condition takes the following form:
\bea
\label{anti0}
\left.\ba{c}
\dis \tilde\dbE\bigg[\l_0\langle\pa_{xx}G(\xi,\mathcal{L}_{\xi})\eta,\eta\rangle+\l_1\langle\pa_{x\mu}G(\xi,\mathcal{L}_{\xi},\tilde\xi)\tilde\eta,\eta\rangle 
\\
\dis+\left|\pa_{xx}G(\xi,\mathcal{L}_{\xi})\eta\right|^2+\l_2\left|\tilde{\mathbb E}_{\mathcal{F}_T^1}[\pa_{x\mu}G(\xi,\mathcal{L}_{\xi},\tilde\xi)\tilde \eta]\right|^2 - \l_3|\eta|^2 \bigg]\leq 0,
\ea\right.
\eea
for some appropriate constants $\l_0>0$, $\l_1\in \dbR$, $\l_2>0$, $\l_3\ge 0$. We remark that the inequality here takes the opposite direction to those in \reff{monsmooth}. In particular, the displacement monotonicity  requires the convexity of $G$ in $x$, while here $G$ is typically concave in $x$, due to the first term in \reff{anti0}. This justifies the name of anti-monotonicity (and to have a better comparison with \reff{monsmooth}, we may also set $\l_1=1$). We also note that, considering the case $\l_3=0$, the second line of \reff{anti0} is positive, this means that the first line of \reff{anti0} should be sufficiently negative, which is exactly in the spirit that $G$ to be sufficiently anti-monotone.  

To establish the global well-posedness of the master equation \reff{master}, we follow the strategy in \cite{GMMZ}, which consists of three steps. The key step of this approach is to show a priori that the anti-monotonicity propagates along the solution $V$. That is, under appropriate conditions, as long as $V(T,\cd)=G$ is anti-monotone, then $V(t,\cd)$ is anti-monotone for all $t$. The second step is to show that the anti-monotonicity of $V$ implies $\pa_x V$ is uniformly Lipschitz continuous in $(x, \mu)$, under $W_2$ in $\mu$. This, together with a representation formula established in \cite{MZ2}, implies further the Lipschitz continuity under $W_1$. In the final step we show that the uniform Lipschitz continuity under $W_1$ enables us to extend a local classical solution to a global one. 

There is a major technical difference from \cite{GMMZ} though. The assumptions we impose for the propagation of anti-monotonicity prevents us from assuming uniform Lipschitz continuity of the data $G$ and $H$. Instead, we can only assume $\pa_x G, \pa_xH$ are uniformly Lipschitz. This has two consequences.  First, the a priori estimate for the boundedness of $\pa_{xx} V$, which is crucial for the global well-posedness of the master equation and is pretty easy to obtain under the conditions in \cite{GMMZ}, becomes very subtle. In fact, we need some serious efforts to obtain this estimate. Moreover, unlike in \cite{GMMZ}, under our conditions the solution $V$ will not be Lipschitz continuous. Instead, we can only expect the Lipschitz continuity of $\pa_x V$. Therefore, we will actually consider the vector master equation of $\vec U:= \pa_x V$ and establish its global well-posedness first. Once we obtain $\vec U$, then it is immediate to solve the original master equation \reff{master} for $V$.

The rest of the paper is organized as follows. In Section \ref{sec:setting} we review the setting in \cite{GMMZ} and introduce our problem. In Section \ref{sec:assum} we introduce the new notion of anti-monotonicity and present the technical conditions used in the paper. In Section \ref{sect-Vantimono} we show a priori the crucial propagation of the  anti-monotonicity. Section \ref{sect-Lipschitz} is devoted to the a priori uniform Lipschitz estimate of $\pa_x V$  in $\mu$, first under $W_2$ and then under $W_1$. In Section \ref{sect-Vxx} we provide the a priori estimate for $\pa_{xx} V$. Finally in Section \ref{sect-wellposedness} we  establish the global well-posedness of the master equation \reff{master}.

\section{The setting} 
\label{sec:setting}
\setcounter{equation}{0}
Throughout the paper we will use the setting in \cite{GMMZ}. We review it briefly in this section and refer to \cite{GMMZ} for more details.

We consider the following product filtered probability space on $[0, T]$: 
\beaa
\O := \O_0 \times \O_1,\q \dbF := \{\cF_t\}_{0\le t\le T} := \{\cF^0_t \otimes \cF^1_t\}_{0\le t\le T},\q \dbP := \dbP_0\otimes \dbP_1,\q \dbE:= \dbE^\dbP.
\eeaa
Here, for $\o = (\o^0, \o^1)\in \O$, $B^0(\o) = B^0(\o^0)$ and $B(\o) = B(\o^1)$ are independent $d$-dimensional Brownian motions; $\dbF^0=\{\cF^0_t\}$ is generated by $B^0$; and $\dbF^1=\{\cF^1_t\}$ is generated by $B$ and $\cF^1_0$, where we assume $\cF^1_0$ has no atom.
Let $(\tilde \O_1, \tilde \dbF^1, \tilde B, \tilde \dbP_1)$ be a copy of the filtered probability space $(\O_1,  \dbF^1, B, \dbP_1)$ and define the larger filtered probability space by
\beaa
\tilde \O := \O\times \tilde\O_1 ,\q \tilde\dbF = \{\tilde \cF_t\}_{0\le t\le T} := \{\cF_t \otimes \tilde \cF^1_t\}_{0\le t\le T},\q \tilde \dbP := \dbP\otimes \tilde\dbP_1,\q \tilde \dbE:= \dbE^{\tilde \dbP}.
\eeaa
Given an $\cF_t$-measurable random variable $\xi(\tilde \o) = \f(\o^0, \o^1)$, $\tilde\o=(\o^0, \o^1, \tilde \o^1)\in \tilde \O$, we see that $\tilde \xi(\tilde \o) := \f(\o^0, \tilde \o^1)$ is a conditionally independent copy of $\xi$, conditional on $\cF^0_t$ under $\tilde \dbP$.  When two conditionally independent copies are needed, we  let $(\bar \O_1, \bar \dbF^1, \bar B, \bar \dbP_1)$ be another copy of $(\O_1,  \dbF^1, B, \dbP_1)$, and enlarge the joint product space further:
\beaa
\bar{\tilde \O}  := \O\times \tilde\O_1\times \bar\O_1,~ \bar{\tilde \dbF} = \{\bar{\tilde \cF}_t\}_{0\le t\le T} := \{\cF_t\otimes\tilde\cF^1_t \otimes \bar \cF^1_t\}_{0\le t\le T},~ \bar {\tilde \dbP} := \dbP\otimes \tilde\dbP_1\otimes \bar\dbP_1, ~\bar{\tilde \dbE}:= \dbE^{\bar{\tilde \dbP}}.
\eeaa
Throughout the paper we will use the probability space $(\O, \dbF, \dbP)$. However, when conditionally independent copies of random variables or processes are needed, we will tacitly use their extensions to the larger space $(\tilde \O, \tilde \dbF, \tilde \dbP)$ $(\bar{\tilde \O}, \bar{\tilde \dbF}, \bar {\tilde \dbP}, \bar{\tilde \dbE})$ without mentioning.

We next introduce the Wasserstein space and differential calculus on Wasserstein space. Let $\cP:=\cP(\dbR^d)$ be the set of all probability measures on $\mathbb R^d$ and, for any $q\ge 1$, let $\cP_q$ denote the set of $\mu\in \cP$ with finite $q$-th moment.
For any sub-$\si$-field $\cG\subset \cF_T$ and $\mu\in \cP_q$, we denote the set of $\dbR^d$-valued, $\cG$-measurable, and $q$-integrable random variables $\xi$ by $\dbL^q(\cG)$;  and the set of $\xi\in \dbL^q(\cG)$ such that the law $\cL_\xi=\mu$ by $\dbL^q(\cG;\mu)$.
For any $\mu,\nu\in \cP_q$,  the $W_q$--Wasserstein distance between them is defined as follows: 
\beaa
W_q(\mu, \nu) := \inf\Big\{\big(\dbE[|\xi-\eta|^q]\big)^{1\over q}:~\mbox{for all $\xi\in \dbL^q(\cF_T; \mu)$, $\eta\in \dbL^q(\cF_T; \nu)$}\Big\}.
\eeaa
For a $W_2$-continuous functions $U: \cP_2 \to \dbR$, its Wasserstein gradient, also called Lions-derivative, takes the form $\pa_\mu U: (\mu,\tilde x)\in \cP_2\times \dbR^d\to \dbR^d$ and satisfies: 
\bea
\label{pamu}
U(\cL_{\xi +  \eta}) - U(\mu) = \dbE\big[\langle \pa_\mu U(\mu, \xi), \eta \rangle \big] + o(\|\eta\|_2), \ \forall\ \xi\in\mathbb L^2(\mathcal{F}_T;\mu),\eta\in\mathbb L^2(\mathcal{F}_T).
\eea
Let $\cC^0(\cP_2)$ denote the set of $W_2$-continuous functions $U:\cP_2\to\dbR$. For $k=1,2$, we introduce $\cC^k(\cP_2)$, which are referred to as functions of {\it full} $\cC^k$ regularity in \cite[Theorem 4.17]{CD1}, as follows. By $\cC^1(\cP_2)$, we mean the space of functions $U\in\cC^0(\cP_2)$ such that $\pa_\mu U$ exists and is continuous on $ \cP_2\times \dbR^d$, it is uniquely determined by \reff{pamu}. Similarly, $\cC^2(\cP_2)$ stands for the set of functions $U\in \cC^1(\cP_2)$ such that $\pa_{\tilde x\mu}U,\pa_{\mu\mu}U$ exist and are continuous on $\cP_2\times\dbR^d$ and $\cP_2\times\dbR^{2d}$ respectively. Let $\cC^2(\dbR^d\times\cP_2)$ denote the set of continuous functions $U:\dbR^d\times\cP_2\to\dbR$ satisfying $\pa_xU,\pa_{xx}U$ exist and are joint continuous on $\dbR^d\times \cP_2$, $\pa_\mu U,\pa_{x\mu}U,\pa_{\tilde x\mu}U$ exist and are continuous on $\dbR^d\times\cP_2\times\dbR^d$, and $\pa_{\mu\mu}U$ exists and is continuous on $\dbR^d\times\cP_2\times\dbR^{2d}$.
%
Finally, we fix the state space 
\[
\Th:= [0, T]\times \dbR^d \times \cP_2
\]
for our master equation, and let $\cC^{1,2}(\Th)$ denote the set of  continuous functions $U\in \Th\to \dbR$ which has the following continuous derivatives: 
$\pa_t U$, $\pa_x U$, $\pa_{xx} U$, $\pa_\mu U$, $  \pa_{x\mu} U$, $\pa_{\tilde x\mu} U,$ $ \pa_{\mu\mu} U.$

One crucial property of  $U\in \cC^{1,2}(\Th)$ functions is  the  It\^{o}  formula. For $i=1,2$, let $d X^i_t := b^i_t dt + \si^i_t dB_t + \si^{i,0}_t d B^0_t,$ where $b^i:[0,T]\times\Omega\to\mathbb R^d$ and $\sigma^i,\sigma^{i,0}:[0,T]\times\Omega\to\mathbb R^{d\times d}$ are $\dbF$-progressively measurable and bounded (for simplicity), and $\rho_t:= \cL_{X^2_t|\cF^0_t}$, then we have
\bea
&&d U(t, X^1_t, \rho_t) =  \Big[\pa_t U + \pa_x U\cdot b^1_t + \frac{1}{2} \tr\big(\pa_{xx} U [\si_t^1 (\si_t^1)^\top + \si_t^{1,0}(\si_t^{1,0})^\top]\big)\Big](t, X^1_t, \rho_t) dt \nonumber\\
&&+\pa_xU(t,X^1_t,\rho_t)\cd\si_t^1dB_t  + (\si^{1,0}_t)^\top\pa_xU(t,X^1_t,\rho_t)\cd dB_t^0
\nonumber\\
&&+\tr \Big(\tilde \dbE_{\cF_t}\big[\pa_\mu U(t,X^1_t,\rho_t,\tilde X^2_t) (\tilde b^{2}_t)^\top\big]\Big) dt +\tilde{\mathbb E}_{\cF_t}\big[(\tilde \si^{2,0}_t)^\top\pa_\mu U(t,X^1_t,\rho_t,\tilde X^2_t) \big]\Big]\cd dB_t^0 \label{Ito}\\
&&+ \tr \Big(\tilde \dbE_{\cF_t}\big[\pa_x\pa_\mu U(t,X^1_t,\rho_t,\tilde X^2_t)\si^{1,0}_t (\tilde \si_t^{2,0})^\top+\frac{1}{2} \pa_{\tilde x}\pa_\mu U(t, X^1_t, \rho_t, \tilde X^2_t)[\tilde \si_t^2 (\tilde \si_t^2)^\top + \tilde \si_t^{2,0}(\tilde \si_t^{2,0})^\top]\big]\nonumber\\
&&\qq +\frac{1}{2}\bar{\tilde \dbE}_{\cF_t}\big[\pa_{\mu\mu}U(t,X^1_t,\rho_t,\tilde X^2_t,\bar X^2_t) \tilde\si_t^{2,0}(\bar \si_t^{2,0})^\top\big]\Big)dt.\nonumber
\eea
See, e.g.,  \cite[Theorem 4.17]{CD2}, \cite{BLPR,CCD}). Here $\cL_{X_t^2|\cF_t^0}$ stands for the conditional law of $X_t^2$ given $\cF_t^0$, and $\tilde \dbE_{\cF_t}$ and $\bar{\tilde \dbE}_{\cF_t}$ are the conditional expectations given $\cF_t$ corresponding to the probability measures $\tilde \dbP$ and $ \bar {\tilde \dbP} $ respectively. Throughout the paper, the elements of $\dbR^d$ are viewed as column vectors; $\pa_x U, \pa_{\mu} U\in\R^d$ are also column vectors; $\pa_{x\mu}U:= \pa_x \pa_\mu U := \pa_x \big[(\pa_\mu U)^\top\big]\in \dbR^{d\times d}$, where $^\top$ denotes the transpose, and similarly for the other second order derivatives; both the notations ``$\cd$'' and $\langle\cdot,\cdot\rangle$ denote the inner product of column vectors. 

We finally introduce the mean field system related to the master equation \reff{master}. It either takes the form of  forward backward McKean-Vlasov SDEs on $[t_0,T]$: given $t_0$ and $\xi\in \dbL^2(\cF_{t_0})$,
\bea
\label{FBSDE1} 
\left.\ba{lll}
\dis\left.\ba{ll}
\dis X^{\xi}_t  \dis =\xi - \int_{t_0}^t\pa_pH(X_s^\xi,\rho_s,Z_s^{\xi})ds+B^{t_0}_t+\b B_t^{0,t_0},\q B^{t_0}_t:= B_t-B_{t_0}, ~B_t^{0,t_0}:= B^0_t-B^0_{t_0}; \\
\dis Y_t^{\xi} =G(X_T^{\xi},\rho_T)+\int_t^T\wh L(X_s^{\xi},\rho_s,Z_s^{\xi})ds- \int_t^TZ_s^{\xi}\cd dB_s-\int_t^TZ_s^{0,\xi}\cd dB_s^{0};\ms
\ea\right.\\
\dis~ \mbox{where}\q \wh L(x,\mu, p) := p\cdot \pa_pH(x,\mu,p)-H(x,\mu,p),\q  \rho_t := \rho^\xi_t:= \cL_{X_t^\xi|\cF^0_t}, 
\ea\right.
\eea
or take the form of forward backward stochastic PDE system on $[t_0, T]$: denoting $\h\beta^2:=1+\b^2$,
\begin{equation}
\label{SPDE}
\left.
\begin{array}{ll}
\dis \!\!\! \!\!\!\!  d \rho(t,x)&\!\!\!\! \dis =  \Big[\frac{\h\beta^2}{2}\tr\big( \pa_{xx} \rho(t,x)\big) + div\big(\rho(t,x) \pa_p H(x, \rho(t,\cd), \pa_x u(t,x))\big)\Big]dt-\b\partial_x\rho(t,x)\cd d B_t^0;\\
\dis \!\!\! \!\!\!\!  d u(t, x)&\!\!\!\! \dis =  v(t,x)\cd dB_t^0 - \Big[\tr\big(\frac{\h\beta^2}{2} \pa_{xx} u(t,x) + \b\partial_x v^\top(t,x)\big) -H(x,\rho(t,\cdot),\pa_x u(t,x))\Big]dt;\\
\dis&\dis \rho(t_0,\cd) = \cL_\xi,\q  u(T,x) = G(x, \rho(T,\cdot)),
\end{array}
\right.
\end{equation}
where the solution triple $(\rho, u, v)$ is $\dbF^0$-progressively measurable and $\rho(t,\cd,\o)$ is a (random) probability measure. 
The systems \reff{FBSDE1} and \reff{SPDE} connect to the master equation \reff{master} as follows: provided all the equations are well-posed and in particular \reff{master} has a classical solution $V$, then 
\bea
\label{YXV}
\left.\ba{c}
Y^\xi_t = V(t, X^\xi_t, \rho_t),\q Z^\xi_t = \pa_x V(t, X^\xi_t, \rho_t),\q \mbox{and}\q u(t,x,\o) = V(t,x, \rho(t,\cd,\o)).
\ea\right.
\eea
It is already well known that, c.f. \cite{CD2}, if the master equation \reff{master} has a classical solution $V$ with bounded derivatives, then we can get existence and uniqueness of the mean field equilibrium, and the equilibrium of the corresponding $N$-player game will converge to the mean field equilibrium. Therefore, we shall only focus on the global well-posedness of the master equation \eqref{master}. 

We conclude this section with the strategy in \cite{GMMZ} for the global well-posedness of \reff{master}. We will follow the same strategy in this paper, except that we shall replace the monotonicity condition with the anti-monotonicity condition:

{\it Step 1.} Introduce appropriate monotonicity condition on data which ensure the propagation of the monotonicity along any classical solution to the master equation. 

{\it Step 2.} Show that the monotonicity of $V(t,\cdot,\cdot)$ implies an (a priori) uniform Lipschitz continuity of $V$ in the measure variable $\mu$.

{\it Step 3.} Combine the local well-posedness of classical solutions and the above uniform Lipschitz continuity to obtain the global well-posedness of classical solutions.

\section{Assumptions and anti-monotonicity conditions} 
\label{sec:assum}
\setcounter{equation}{0}

In this section, we introduce the following notations. For any $A\in \dbR^{d\times d}$,
\bea
\label{kappaA}
\left.\ba{c}
\dis \underline \k(A) := \inf_{|x|=1} \langle Ax, x\rangle = \mbox{the smallest eigenvalue of ${1\over 2}[A+A^\top]$},\qq \overline \k(A) :=  \sup_{|x|=1} \langle Ax, x\rangle;\\
\dis \underline \k'(A) := \text{ the smallest real part of eigenvalues of $A$};\\
\dis |A| :=  \sup_{|x|=|y|=1} \langle Ax, y\rangle.
\ea\right.
\eea
It is obvious that, for any $A,A_1,A_2\in \mathbb R^{d\times d}$ and $x\in\mathbb R^d$,
\bea
\label{kappaA2}
\left.\ba{c}
\dis \mbox{$|\cd|$ is a norm on $\dbR^{d\times d}$, \q $|A_1A_2|\le |A_1||A_2|$,\q $|Ax|\le |A||x|$},\\
\dis \mbox{and, when $A$ is symmetric, $\underline \k'(A) = \underline \k(A)$ and $|A|= |\underline\k(A)|\vee |\overline\k(A)|$}.
\ea\right.
 \eea

\subsection{Regularity assumptions}

We first specify some technical assumptions on $G$ and $H$.

\begin{assum} \label{assum-regH} 
(i) $H\in \cC^2(\mathbb R^d\times\mathcal{P}_2\times\mathbb R^d)$ and there exist constants $\overline L_{xp}^{H},\overline L_{xx}^{H}$,$L_2^{H}>0$ such that 
\begin{equation}
\label{Hbound}
|\pa_{xp}H|\leq \overline L_{xp}^{H},\quad |\pa_{xx}H|\leq \overline L^{H}_{xx},\qq
|\pa_{pp}H|,|\pa_{x\mu}H|,|\pa_{p\mu}H|\le L_2^{H}.
\end{equation}
(ii) $H\in \cC^3(\mathbb R^d\times\mathcal{P}_2\times\mathbb R^d)$, and
\beaa
&\dis \pa_xH,\;  \pa_pH,\; \pa_{xx}H,\;  \pa_{xp}H,\;  \pa_{pp}H,\;  \pa_{xxp}H,\;  \pa_{xpp}H,\;  \pa_{ppp}H\; \in \cC^2(\mathbb R^d\times\mathcal{P}_2\times\mathbb R^d), \\
&\dis \pa_\mu H,\pa_{x\mu}H, \pa_{p\mu}H,\pa_{ xp\mu}H,\pa_{ pp\mu}H\in \cC^2(\mathbb R^d\times\mathcal{P}_2\times\mathbb R^{2d}),
\eeaa
 where all the second and higher order derivatives of $H$ involved above are uniformly bounded. 
 \end{assum}

\begin{assum}\label{assum-regG}
(i) $G\in \cC^2(\mathbb R^d\times\mathcal{P}_2)$, and there exist  constants $\overline L_{xx}^G, L_2^G>0$ such that 
\bea
\label{Gbound}
| \pa_{xx}G|\le \overline L^G_{xx},\q |\pa_{x\mu}G|\leq L_2^G.
\eea
(ii) 
$
\pa_xG, \pa_{xx}G\in \cC^2(\mathbb R^d\times\mathcal{P}_2)$, and $\pa_\mu G, \pa_{x\mu}G \in \cC^2(\mathbb R^d\times\mathcal{P}_2\times\mathbb R^d),
$
and all the second and higher order derivatives of $G$ involved here are uniformly bounded. 
\end{assum}

Here the spaces $\cC^2, \cC^3$ are defined in the same manner as $\cC^{1,2}(\Th)$.  Note that at above we do not require the first order derivatives to be uniformly bounded. In fact, the condition \reff{antiH1} below does not allow $\pa_x H$ to be bounded.

\begin{rem}\label{rem:W2Lip}
Under Assumption \ref{assum-regG}-(i), we see that $\pa_xG$ is uniformly Lipschitz continuous in $\mu$ under $W_1$ on $\mathbb R^d\times\mathcal{P}_2$ with Lipschitz constant $L_2^G$. This implies further the Lipschitz continuity of $\pa_x G$ in $\mu$ under $W_2$ on $\mathbb R^d\times\mathcal{P}_2$, and we denote the Lipschitz constant  by $\tilde L_2^G\leq L_2^G$:
\beaa
\tilde {\mathbb E}\Big[[\pa_{x\mu}G(x,\mu,\tilde \xi)\tilde \eta]\Big]\leq \tilde L_2^G\Big(\mathbb E[|\eta|^2]\Big)^{\frac{1}{2}},\qq \forall \xi\in\mathbb L^2(\mathcal{F}_T^1,\mu),\eta\in\mathbb L^2(\mathcal{F}_T^1).
\eeaa
\end{rem}

\subsection{Monotonicity and anti-monotoncity conditions}

Under the above regularity conditions on the data $G$ and $H$, the MFG may still have multiple mean field equilibria  over a long time duration and thus the global well-posedness of classical solutions for the master equations can fail. Therefore, some  structural conditions on $G, H$ are needed in order to guarantee its global well-posedness. The typical structural conditions assumed in the literature are two types of monotonicity conditions, i.e., the Lasry-Lions monotonicity condition and the displacement monotonicity condition. 

\begin{defn} 
\label{defn-mon}
Let $U :\dbR^d\times \cP_2\to\mathbb R$ be such that $U\in \cC^2(\dbR^d\times \cP_2)$.
\begin{enumerate}
\item[(i)] $U$ is called \emph{Lasry-Lions monotone}, if  for any $\xi, \eta\in \dbL^2(\mathcal{F}^1_T)$,
\bea
\label{mon1}
\tilde{\mathbb E}\Big[\langle\pa_{x\mu}U(\xi,\mathcal{L}_{\xi},\tilde\xi)\tilde\eta,\eta\rangle\Big] \ge 0.
\eea

\item[(ii)] $U$ is called \emph{displacement monotone} if for any $\xi, \eta\in \dbL^2(\mathcal{F}^1_T)$,
\bea
\label{displacement1}
\tilde\dbE\Big[\langle\pa_{x\mu}U(\xi,\mathcal{L}_{\xi},\tilde \xi)\tilde\eta,\eta\rangle+\langle\pa_{xx}U(\xi,\mathcal{L}_{\xi})\eta,\eta\rangle\Big] \ge 0.
\eea
\item[(iii)] $U$ is called \emph{displacement semi-monotone} if, for some $\l\in\dbR$ and for any $\xi, \eta\in \dbL^2(\mathcal{F}^1_T)$,
\bea
\label{semi-displacement}
\tilde\dbE\Big[\langle\pa_{x\mu}U(\xi,\mathcal{L}_{\xi},\tilde \xi)\tilde\eta,\eta\rangle+\langle\pa_{xx}U(\xi,\mathcal{L}_{\xi})\eta,\eta\rangle\Big] - \l \dbE[|\eta|^2] \ge 0.
\eea
\end{enumerate}
\end{defn}
Here, as in Section \ref{sec:setting}, $(\tilde\xi,\tilde\eta)$ is an independent copy of $(\xi,\eta)$. We remark that the displacement semi-monotonicity is obviously weaker than the  displacement monotonicity \reff{displacement1}, and when $\pa_{xx} U$ is bounded, it is also weaker than the Lasry-Lions monotonicity \reff{mon1}.

\begin{rem}\label{rem:implicationDisp} The above formulations of the monotonicity conditions are convenient for our purpose. For $U \in \cC^2(\dbR^d\times \cP_2)$, \reff{mon1} and \reff{displacement1} are equivalent to \reff{LLmon} and \reff{dismon}, respectively, which appear more often in the literature. See \cite[Remark 2.4]{GMMZ}. 
 \end{rem}

We next turn to the monotonicity conditions for the Hamiltonian $H$. In the literature, the Lasry-Lions monotonicity has only been proposed for the separable Hamiltonians, i.e., $H(x, \mu, p) = H_0(x, p)-F(x, \mu)$ and $F$ satisfies \eqref{LLmon}.
In \cite{GMMZ}, a notion of displacement monotonicity for non-separable $H$ was proposed to study the well-posedness of the master equation \eqref{master}. 
\begin{defn} \label{defn:displace-mono-H} Let $H$ be a Hamiltonian satisfying Assumption \ref{assum-regH}(i) and $H$ is strictly convex in $p$. We say that $H$ is displacement monotone if: for any $\xi,\eta\in \dbL^2(\cF^1_T)$ and any bounded Lipschitz continuous function $\f\in C^1(\dbR^d; \dbR^d)$, 
\bea\label{Hconvex} 
\left.\ba{c}
\dis \tilde \dbE\bigg[\Big\langle \pa_{x\mu} H(\xi, \cL_{\xi}, \tilde \xi, \f(\xi))\tilde \eta + \pa_{xx} H(\xi, \cL_{\xi}, \f(\xi)) \eta,~\eta\Big\rangle \\
\dis+{1\over 4}\Big| \big(\pa_{pp} H(\xi, \cL_{\xi}, \f(\xi))\big)^{-{1\over 2}}\tilde \dbE_{\cF^1_T}[\pa_{p\mu} H(\xi, \cL_{\xi}, \tilde \xi, \f(\xi)) \tilde \eta] \Big|^2\bigg]\le 0.
\ea\right.
\eea 
\end{defn}
\begin{rem}
(i) The above definition of displacement monotonicity for non-separable Hamiltonians is not really used in the rest of the paper except for the comparison with the new notion of anti-monotonicity introduced below. We refer to \cite[Proposition 3.7]{GMMZ} for another equivalent definition of the above one.

(ii) The function $\varphi(\xi)$ in \eqref{Hconvex} is chosen to be $\pa_xV(t,\xi,\cL_{\xi})$ in the proof of the propagation of the displacement monotonicity \reff{displacement1} along $V(t,\cd)$ in \cite{GMMZ}. Since $\pa_xV$ is not priorily known, the displacement monotonicity \eqref{Hconvex} is made for any desirable function $\varphi$.

(iii) When $H$ is non-separable,  it still remains a challenge to find appropriate conditions on $H$ so that the Lasry-Lions monotonicity \reff{mon1} could propagate along the solution $V(t,\cd)$.
\end{rem}

Finally we introduce the anti-monotonicity condition, which is the main structural condition in this paper and serves as an alternative sufficient  condition for the global wellposedness of the master equation.  Denote 
\bea
\label{D4}
D_4 := \Big\{\vec \l = (\l_0, \l_1, \l_2, \l_3): \l_0>0, \l_1\in \dbR, \l_2 >0, \l_3\ge 0\Big\}.
\eea

\begin{defn}
\label{defn-anti}
Let $U\in \cC^2(\dbR^d\times \cP_2)$ and $\vec\l\in D_4$. We say $U$ is \mbox{$\vec\l$-anti-monotone} if,
\bea
\label{anti}
\left.\ba{c}
(AntiMon)^{\vec \l}_\xi U(\eta, \eta):=\tilde\dbE\bigg[\l_0\langle\pa_{xx}U(\xi,\mathcal{L}_{\xi})\eta,\eta\rangle+\l_1\langle\pa_{x\mu}U(\xi,\mathcal{L}_{\xi},\tilde\xi)\tilde\eta,\eta\rangle \\
+\left|\pa_{xx}U(\xi,\mathcal{L}_{\xi})\eta\right|^2+\l_2\left|\tilde{\mathbb E}_{\mathcal{F}_T^1}[\pa_{x\mu}U(\xi,\mathcal{L}_{\xi},\tilde\xi)\tilde \eta]\right|^2 - \l_3|\eta|^2 \bigg]\leq 0,\qq \forall  \xi, \eta\in \dbL^2(\mathcal{F}^1_T).
\ea\right.
\eea
\end{defn}

\begin{rem}
\label{rem-anti}
(i) The main feature of \reff{anti} is that the direction of the inequality is opposite to those in Definition \ref{defn-mon}. In particular, \eqref{anti} implies the Lasry-Lions anti-monotonicity, i.e.
\begin{equation}\label{eq:LLantimonotonicity}
\tilde\dbE\Big[\langle\pa_{x\mu}U(\xi,\mathcal{L}_{\xi},\tilde\xi)\tilde\eta,\eta\rangle \Big]\leq 0,
\end{equation}
for the case that $\l_0=\l_3=0$ and $\l_1=\l_2 =1$. In fact, in this case the condition \eqref{anti} is  stronger than \eqref{eq:LLantimonotonicity} and we interpret it as $U$ is sufficiently Lasry-Lions anti-monotone:
\begin{equation}\label{eq:LLantimonotonicity2}
\tilde\dbE\Big[\langle\pa_{x\mu}U(\xi,\mathcal{L}_{\xi},\tilde\xi)\tilde\eta,\eta\rangle \Big] \le -  \tilde\dbE\bigg[\left|\pa_{xx}U(\xi,\mathcal{L}_{\xi})\eta\right|^2+\left|\tilde{\mathbb E}_{\mathcal{F}_T^1}[\pa_{x\mu}U(\xi,\mathcal{L}_{\xi},\tilde\xi)\tilde \eta]\right|^2  \bigg]\leq 0.
\end{equation}

Similarly, in the case $\l_0=\l_1=\l_2=1$ and $\l_3=0$, we see that  \eqref{anti} implies $U$ is sufficiently displacement anti-monotone:
{\small \begin{equation}\label{eq:antimonotonicity}
\tilde\dbE\Big[\langle\pa_{xx}U(\xi,\mathcal{L}_{\xi})\eta,\eta\rangle+\langle\pa_{x\mu}U(\xi,\mathcal{L}_{\xi},\tilde\xi)\tilde\eta,\eta\rangle \Big]  \le -  \tilde\dbE\bigg[\left|\pa_{xx}U(\xi,\mathcal{L}_{\xi})\eta\right|^2+\left|\tilde{\mathbb E}_{\mathcal{F}_T^1}[\pa_{x\mu}U(\xi,\mathcal{L}_{\xi},\tilde\xi)\tilde \eta]\right|^2  \bigg]\leq 0.
\end{equation}}
 Note that the concavity of $U$ in $x$ could help in \eqref{eq:antimonotonicity}, while in \reff{displacement1} its convexity is helpful.\\

(ii) The inequality \eqref{anti} implies the displacement semi-anti-monotonicity, i.e.
\begin{equation}\label{eq:semiantimon}
\tilde\dbE\Big[\langle\pa_{xx}U(\xi,\mathcal{L}_{\xi})\eta,\eta\rangle+\langle\pa_{x\mu}U(\xi,\mathcal{L}_{\xi},\tilde\xi)\tilde\eta,\eta\rangle \Big]\leq \l_3\tilde\dbE\big[|\eta|^2\big],
\end{equation}
if $\vec{\l}\in D_4$, $\l_0=\l_1=1$ and $\l_3\geq 0$. Note that the condition \eqref{eq:semiantimon} is weaker than \eqref{eq:antimonotonicity} for the case. We recall that in the literature a function $u:\dbR^d\to\dbR$ is said to be semi-concave, or $\l$-concave, if  $\pa_{xx}u\leq \l I_d $ for some constant $\l>0$, 
where $I_d$ is the $d\times d$ identity matrix. 
We follow the same spirit to call $U$  \mbox{$\vec\l$-anti-monotone} if $U$ satisfies \eqref{anti}.
\end{rem}

We next provide an example which is $\vec\l$-anti-monotone.
\begin{eg}
\label{eg-anti}
Let $d=1$ and consider the function: for some constants $a_0, a_1$, 
\[
 U(x,\mu)= {a_0\over 2}|x|^2 + a_1 x\int_{\dbR}y\mu(dy),\q (x,\mu)\in \dbR \times \cP_2.
\]
It is clear that  $\pa_{xx}U= a_0$ and $\pa_{x\mu} U= a_1$. 

(i) For any $\xi, \eta\in \dbL^2(\mathcal{F}^1_T)$, we have
\[
\tilde{\mathbb E}\Big[\langle\pa_{x\mu} U(\xi,\mathcal{L}_{\xi},\tilde\xi)\tilde\eta,\eta\rangle\Big]= a_1\big|\dbE[\eta]\big|^2.
\]
So $U$ is  Lasry-Lions monotone  if $a_1\ge 0$, and Lasry-Lions anti-monotone if $a_1 \le 0$.

(ii) Similarly we have
\[
\tilde\dbE\Big[\langle\pa_{xx} U(\xi,\mathcal{L}_{\xi})\eta,\eta\rangle+\langle\pa_{x\mu} U(\xi,\mathcal{L}_{\xi},\tilde\xi)\tilde\eta,\eta\rangle \Big]=a_0\dbE\big[|\eta|^2\big] + a_1\big|\dbE[\eta]\big|^2.
\]
Then one can easily check that $U$ is displacement monotone if $a_0 \ge 0$,  $a_1 \ge -a_0$, and displacement anti-monotone if $a_0\le 0$, $a_1 \le -a_0$. 

(iii) For any $\vec\l\in D_4$, we have
\beaa
(AntiMon)^{\vec \l}_\xi U(\eta, \eta):=\Big[\l_0a_0 + |a_0|^2 - \l_3\Big]\dbE[|\eta|^2] +\big[\l_1a_1 + \l_2|a_1|^2\big] \big|\dbE[\eta]\big|^2.
\eeaa
Then $U$ is $\vec\l$-anti-monotone if:
\beaa
\l_0a_0 + |a_0|^2 - \l_3 \le 0,\q \l_0a_0 + |a_0|^2  - \l_3 \le -\big[\l_1a_1 + \l_2|a_1|^2\big],
\eeaa
which is equivalent to:
\beaa
\l_3 \ge \max\Big(\l_0a_0 + |a_0|^2,~~ \l_0a_0 + |a_0|^2+\l_1a_1 + \l_2|a_1|^2\Big).
\eeaa
 In particular, if we set $\l_0 = \l_1 = \l_2 = 1$, $\l_3=0$, and $-1 \le a_0, a_1 \le 0$, we see that $U$ is $\vec\l$-anti-monotone.
%
\end{eg}
\begin{rem}
Let $U\in \cC^2(\dbR^d\times\cP_2)$ and $\vec{\l}\in D_4$.

(i) When $\l_0=0$, \eqref{anti} is equivalent to the following integral form: for any $\xi_1,\xi_2\in \dbL^2(\mathcal{F}^1_T)$,
\begin{eqnarray}
&& \l_1\dbE\Big[ U(\xi_1, \cL_{\xi_1}) + U(\xi_2, \cL_{\xi_2}) -U(\xi_1,\cL_{\xi_2})- U(\xi_2,\cL_{\xi_1})\Big]\nonumber\\
&+&\dbE\Big[|\pa_{x}U(\xi_1,\cL_{\xi_2})-\pa_xU(\xi_2,\cL_{\xi_2})|^2+\l_2|\pa_xU(\xi_2,\cL_{\xi_1})-\pa_xU(\xi_2,\cL_{\xi_2})|^2\Big]\\
&\leq& \l_3 \dbE\big[|\xi_1-\xi_2|^2\big]+o\big(\dbE\big[|\xi_1-\xi_2|^2\big]\big).\nonumber
\end{eqnarray}
Here $o(\e)$ means it vanishes faster than $\e$ as $\e\to 0$.

(ii) When $\l_0=\l_1$,  \eqref{anti} is equivalent to the following integral form: for any $\xi_1,\xi_2\in \dbL^2(\mathcal{F}^1_T)$,
\begin{eqnarray}
&&\l_0 \dbE\Big[\big\langle\pa_xU(\xi_1,\cL_{\xi_1})-\pa_xU(\xi_2,\cL_{\xi_2}),\xi_1-\xi_2\big\rangle\Big]\nonumber\\
&+&\dbE\Big[|\pa_{x}U(\xi_1,\cL_{\xi_2})-\pa_xU(\xi_2,\cL_{\xi_2})|^2+\l_2|\pa_xU(\xi_2,\cL_{\xi_1})-\pa_xU(\xi_2,\cL_{\xi_2})|^2\Big]\\
&\leq&C\dbE\big[|\xi_1-\xi_2|^2\big]+o\big(\dbE\big[|\xi_1-\xi_2|^2\big]\big).\nonumber
\end{eqnarray}

(iii) In general, \eqref{anti} is equivalent to the following integral form: for any $\xi_1,\xi_2\in \dbL^2(\mathcal{F}^1_T)$,
\begin{eqnarray}
&&\dbE\Big[\l_0\big\langle\pa_xU(\xi_1,\cL_{\xi_2})-\pa_xU(\xi_2,\cL_{\xi_2}),\xi_1-\xi_2\big\rangle+\l_1\big\langle\pa_xU(\xi_2,\cL_{\xi_1})-\pa_xU(\xi_2,\cL_{\xi_2}),\xi_1-\xi_2\big\rangle\Big]\nonumber\\
&+&\dbE\Big[|\pa_{x}U(\xi_1,\cL_{\xi_2})-\pa_xU(\xi_2,\cL_{\xi_2})|^2+\l_2|\pa_xU(\xi_2,\cL_{\xi_1})-\pa_xU(\xi_2,\cL_{\xi_2})|^2\Big]\\
&\leq& \l_3 \dbE\big[|\xi_1-\xi_2|^2\big]+o\big(\dbE\big[|\xi_1-\xi_2|^2\big]\big).\nonumber
\end{eqnarray}

\end{rem}

\begin{assum}\label{assum-antidisplacement} 
(i) $G$ satisfies Assumption \ref{assum-regG}-(i) and is $\vec \l$-anti-monotone for some  $\vec\l\in D_4$;

(ii) $H$ satisfies Assumption \ref{assum-regH}-(i) and there exist constants $\underline L_{xp}^{H}>0,\underline L_{xx}^{H}>0$, $\overline \g >\underline \g  >0$ s.t.
\bea\label{antiH1}
&\dis \underline\k(\pa_{xp}H)\geq \underline L_{xp}^{H},\q  \underline\k(\pa_{xx}H)\geq\underline{L}_{xx}^{H},\\
\label{antiH2}
&\dis \underline \g \underline L^H_{xp} \le \underline L^H_{xx} \le \overline L^H_{xx} \le \overline \g \underline L^H_{xp},\q\overline L_{xp}^{H} \le \overline \g \underline L^H_{xp}.
\eea
\end{assum}
Note that we do not require structural conditions on $\pa_{x\mu} H$ here, and $\pa_{pp} H$ can be degenerate.

\section{Propagation of anti-monotonicity}
\label{sect-Vantimono}
\setcounter{equation}{0}

In this section we show that any classical solution $V$ to the master equation \eqref{master} could propagate the anti-monotonicity  under appropriate conditions. 

\begin{thm}\label{thm-antidisplacement}
Let Assumption \ref{assum-antidisplacement} hold and $V$ be a classical solution of the master equation \eqref{master} such that 
$$\pa_{xx}V(t,\cdot,\cdot)\in \cC^2(\mathbb R^d\times\mathcal{P}_2), \quad  \pa_{x\mu}V(t,\cdot,\cdot,\cdot)\in \cC^2(\mathbb R^d\times\mathcal{P}_2\times\mathbb R^d),$$ 
and all the second and higher order derivatives of $V$ involved above are also continuous in the time variable and are uniformly bounded. Assume further that there exist a constant $L^V_{xx} >0$  such that
\bea
\label{Vxx}
&\dis |\pa_{xx} V|\le L^V_{xx},\\
\label{theta}
&\dis\text{and}\quad \l_0 > { \overline \g^2 [1+L^V_{xx}]^2-8\l_3  \over 4 \underline \g} \q\mbox{so that}\q  \th_1 := { \overline \g [1+L^V_{xx}]  \over \sqrt{4( \underline \g \l_0+2\l_3)}} <1.
\eea
Introduce the following symmetric matrices, which depend only on $\underline \g, \overline \g, \vec \l$, and $L^V_{xx}$:
\bea
\label{A}
\left.\ba{c}
A_1:= 
\begin{bmatrix}
4[1-\th_1]   & 0& 0 \\
0 &  2\l_2 & 0\\
0&  0& [1-\th_1][\l_0\underline \g+2\l_3] 
\end{bmatrix},\ms\\
A_2:= \begin{bmatrix}
 \l_0 & \l_0 & |\l_0-{1\over 2}\l_1|+\l_3 \\
\l_0  &  | \l_1| &{1\over 2}|\l_1|+ \l_2 +\l_3 \\
|\l_0-{1\over 2}\l_1|+\l_3  & {1\over 2}|\l_1|+\l_2 + \l_3 &  |\l_1|+2\l_3
\end{bmatrix}  +  \begin{bmatrix}
0 & 1 & 1\\
1 &  \l_2 &\l_2\\
1 &  \l_2& 0
\end{bmatrix} L_{xx}^V.
\ea\right.
\eea
 Then, whenever 
 \bea
 \label{Lxp}
 \underline L^H_{xp} \ge \underline \k(A_1^{-1}A_2) L^H_2,
 \eea
  $V(t,\cd)$ is $\vec \l$-anti-monotone in the sense of \reff{anti} for all $t\in [0, T]$. 
\end{thm}
\proof Without loss of generality, we shall prove the theorem only for $t_0=0$. 

Fix $\xi\in \dbL^{2}(\cF_0)$ and $\eta\in \dbL^{2}(\cF_0)$. Given the desired regularity of $V$ and $H$, the following system of McKean-Vlasov SDEs has a unique solution $(X, \d X)$: 
\bea\label{XY}
\left.\ba{lll}
\dis X_t = \xi -\int_0^t \pa_pH(X_s, \mu_s, \pa_x V(s, X_s, \mu_s)) ds +  B_t+\beta B_t^0,\q \mu_t := \cL_{X_t|\mathcal{F}_t^{0}};\\
\dis \delta X_t = \eta -\int_0^t \Big[\pa_{px}H (X_s,\mu_s,\pa_xV(s,X_s,\mu_s)) \delta  X_s + \tilde \dbE_{\mathcal{F}_s}[ \pa_{p\mu}H(X_s,\mu_s,\pa_xV(s,X_s,\mu_s),\tilde X_s) \delta  \tilde X_s]\\
\dis\qq\qq\qq\q+ \pa_{pp}H(X_s,\mu_s,\pa_xV(s,X_s,\mu_s)) [\Upsilon_s+\bar\Upsilon_s]\Big]ds,\\
\dis \mbox{where}\q \Upsilon_t :=\tilde \dbE_{\mathcal{F}_t}[\pa_{x\mu}  V(t,X_t,\mu_t, \tilde X_t) \delta \tilde X_t], \q \bar\Upsilon_t := \pa_{xx} V(t,X_t,\mu_t)  \delta  X_t.
\ea\right.
\eea
In the sequel,  for simplicity of notation, we omit the variables $(t, \mu_t)$ as well as the dependence on $\partial_x V$, and denote 
\beaa
H_p(X_t) := \pa_p H\big(X_t, \mu_t, \pa_x V(t, X_t, \mu_t)\big), \q H_{p\mu}(X_t, \tilde X_t) := \pa_{p\mu} H\big(X_t, \mu_t, \tilde X_t, \pa_x V(t, X_t, \mu_t)\big),
\eeaa
and similarly for $H_{xp}, H_{pp}$,  $H_{x\mu}$, $\partial_{xx}V$, $\partial_{x\mu}V$, etc. 
We remark that, $(\tilde X_t, \d\tilde X_t)$ is a conditionally independent copy of $(X_t, \d X_t)$ and $\mu_t$ is $\cF^0_t$-measurable.

Recall \reff{XY} and introduce: 
\bea
\label{IUpsilon}
\left.\ba{c}
\dis I_t:=\dbE\big[\langle \Upsilon_t, \delta  X_t\rangle \big], \q \bar I_t:=\dbE\big[\langle \bar \Upsilon_t, \delta  X_t\rangle \big];\ms\\
\G_t :=  (AntiMon)^{\vec \l}_{X_t}V(t,\cd)(\d X_t, \d X_t) =\l_0 \bar I_t+ \l_1 I_t+   \mathbb E\big[|\bar \Upsilon_t|^2+\l_2 |\Upsilon_t|^2-\l_3 |\delta X_t|^2\big].
\ea\right.
\eea
By the calculation in \cite[Theorem 4.1]{GMMZ} we have
\bea
\label{Idot}
\left.\ba{lll}
\dis {d\over dt} I(t) = \dbE\Big[-\big\langle H_{pp}(X_t) \Upsilon_t,~\Upsilon_t\big\rangle -\big\langle \tilde \dbE_{\mathcal{F}_t}\big[ H_{p\mu}(X_t,\tilde X_t)\delta  \tilde X_t\big] , \Upsilon_t - \bar \Upsilon_t\big\rangle  \\
\dis\qq\qq + \big\langle \tilde{\mathbb E}_{\mathcal{F}_t}\big[H_{x\mu}(X_t,\tilde X_t) \delta  \tilde X_t\big], \delta  X_t\big\rangle\Big];\\
\dis{d\over dt} {\bar I}(t) = \dbE\Big[-\big\langle H_{pp}(X_t)\bar \Upsilon_t,\bar \Upsilon_t\big\rangle-2\big\langle H_{pp}(X_t)\bar \Upsilon_t, \Upsilon_t\big\rangle\\
\dis\qq\qq - 2\big\langle\bar \Upsilon_t,\tilde\dbE_{\mathcal{F}_t}[H_{p\mu}(X_t,\tilde X_t)\delta \tilde X_t]\big\rangle+\big\langle H_{xx}(X_t)\delta X_t,\delta X_t\big\rangle \Big],
\ea\right.
\eea
and, by the calculation in \cite[Theorem 5.1]{GMMZ} we have
\bea\label{dUpsilon}
\left.\ba{lll}
\dis d\Upsilon_t =(dB_t)^\top K_1(t) + \b  (dB^0_t)^\top K_2(t)+ \big[ K_3(t) \Upsilon_t + K_4(t)\big] dt   ;\\
d\bar\Upsilon_t = (dB_t)^\top \bar K_1(t)+ \b  (dB^0_t)^\top \bar K_2(t)+ \big[ 2H_{xp}(X_t) \bar\Upsilon_t-\pa_{xx}V(X_t)H_{pp}(X_t)\Upsilon_t + \bar K_3(t)\big] dt  ,
\ea\right.
\eea
where ($K_5(t)$ and $K_6(t)$ in \cite{GMMZ} turn to $K_3(t)$ and $K_4(t)$ respectively here)
\bea
\label{K}
\left.\ba{lll}
K_1(t) := \tilde{\mathbb E}_{\mathcal{F}_t}\big[\pa_{xx\mu}V(X_t,\tilde X_t)\delta \tilde X_t\big],\\
 K_2(t):= K_1(t) + \bar{\tilde \dbE}_{\mathcal{F}_t}\Big[\big[(\pa_{\mu x\mu}V)(X_t,\bar X_t,\tilde X_t)  +\pa_{\tilde xx\mu}V(X_t,\tilde X_t)\big]\delta\tilde X_t\Big],\\
K_3(t) := H_{xp}(X_t) + \pa_{xx} V(X_t) H_{pp}(X_t),\\
  K_4(t) := \tilde\dbE_{\mathcal{F}_t}\Big[  \big[H_{x\mu}(X_t,\tilde X_t) +  \pa_{xx} V(X_t)H_{p\mu}(X_t,\tilde X_t) \big]\delta\tilde X_t\Big], \\
   \bar K_1(t) := \pa_{xxx}V(X_t)\delta  X_t,\\
 \bar K_2(t):= \bar K_1(t) +  \tilde \dbE_{\mathcal{F}_t}\Big[(\pa_{\mu xx}V)(X_t,\tilde X_t)\delta\tilde X_t\Big],\\
 \bar K_3(t) :=[H_{xx}(X_t)-\pa_{xx}V(X_t)H_{px}(X_t)]\delta X_t- \pa_{xx}V(X_t)\tilde\dbE_{\mathcal{F}_t} \big[H_{p\mu}(X_t,\tilde X_t)\delta \tilde X_t\big].\\
\ea\right.
\eea
In particular, this implies that
\bea\label{dUpsilon2}
\left.\ba{lll}
\dis {d\over dt} \dbE[|\Upsilon_t|^2]  \ge 2\dbE\Big[\big\langle \Upsilon_t,~  K_3(t) \Upsilon_t + K_4(t)\big\rangle\Big];\ss\\
\dis {d\over dt} \dbE[|\bar\Upsilon_t|^2] \ge 2 \dbE\Big[\big\langle \bar\Upsilon_t ,~ 2H_{xp}(X_t) \bar\Upsilon_t-\pa_{xx}V(X_t)H_{pp}(X_t)\Upsilon_t + \bar K_3(t)\big\rangle\Big].
\ea\right.
\eea
Moreover,  by  \eqref{XY} we have
\begin{equation}\label{deltaX2}
\frac{d}{dt}\mathbb E\left[|\delta X_t|^2\right]=-2\mathbb E\Big[\Big\langle H_{px}(X_t)\delta X_t+\tilde{\mathbb E}_{\mathcal{F}_t}[H_{p\mu}(X_t,\tilde{X_t})\delta \tilde X_t]+H_{pp}(X_t)[ \Upsilon_t+\bar \Upsilon_t],~\delta X_t \Big\rangle\Big].
\end{equation}
Thus, by \reff{Idot}, \reff{dUpsilon2}, and \reff{deltaX2},  we have
\beaa
\left.\ba{lll}
\dis  {d\over dt}\G_t  \ge \l_0 \dbE\Big[-\big\langle H_{pp}(X_t)\bar \Upsilon_t,\bar \Upsilon_t\big\rangle-2\big\langle H_{pp}(X_t)\bar \Upsilon_t, \Upsilon_t\big\rangle\\
\dis\qq\qq - 2\big\langle\bar \Upsilon_t,\tilde\dbE_{\mathcal{F}_t}[H_{p\mu}(X_t,\tilde X_t)\delta \tilde X_t]\big\rangle+\big\langle H_{xx}(X_t)\delta X_t,\delta X_t\big\rangle \Big] \\
\dis \qq + \l_1  \dbE\Big[-\big\langle H_{pp}(X_t) \Upsilon_t,~\Upsilon_t\big\rangle -\big\langle \tilde \dbE_{\mathcal{F}_t}\big[ H_{p\mu}(X_t,\tilde X_t)\delta  \tilde X_t\big] , \Upsilon_t - \bar \Upsilon_t\big\rangle  \\
\dis\qq\qq + \big\langle \tilde{\mathbb E}_{\mathcal{F}_t}\big[H_{x\mu}(X_t,\tilde X_t) \delta  \tilde X_t\big], \delta  X_t\big\rangle\Big]\\
\dis\qq +2\dbE\Big[ \big\langle \bar\Upsilon_t, \big[ 2H_{xp}(X_t) \bar\Upsilon_t-\pa_{xx}V(X_t)H_{pp}(X_t)\Upsilon_t + \bar K_3(t)\big] \big\rangle +\l_2\big\langle\Upsilon_t, \big[ K_3(t) \Upsilon_t + K_4(t)\big] \big\rangle\Big]\\
\dis\qq +2\l_3\mathbb E\Big[\big\langle H_{px}(X_t)\delta X_t+\tilde{\mathbb E}_{\mathcal{F}_t}[H_{p\mu}(X_t,\tilde{X_t})\delta \tilde X_t]+H_{pp}(X_t)[ \Upsilon_t+\bar \Upsilon_t],\delta X_t \big\rangle\Big]\\
\dis = \dbE\Big[\big\langle [ -\l_0 H_{pp}(X_t)+4 H_{xp}(X_t)]  \bar \Upsilon_t, \bar \Upsilon_t \big\rangle +\big \langle[ - \l_1H_{pp}(X_t) + 2\l_2K_3(t)] \Upsilon_t, \Upsilon_t\big\rangle \\
\dis\qq +\big\langle  [\l_0 H_{xx}(X_t)+2\l_3H_{px}(X_t)] \d X_t,\d X_t\big\rangle\\
\dis \qq +\big\langle  \l_1\tilde{\mathbb E}_{\mathcal{F}_t}\big[H_{x\mu}(X_t,\tilde X_t)  \delta  \tilde X_t\big]+2\l_3\tilde{\mathbb E}_{\mathcal{F}_t}\big[H_{p\mu}(X_t,\tilde X_t)  \delta  \tilde X_t\big],~\d X_t\big\rangle\\
\dis\qq -\big\langle  2 [\l_0 H_{pp}(X_t) + \pa_{xx} V(X_t) H_{pp}(X_t)]\Upsilon_t,\bar \Upsilon_t\big\rangle \\
\dis\qq   +\big\langle [-2\l_0+\l_1] \tilde \dbE_{\mathcal{F}_t}\big[ H_{p\mu}(X_t,\tilde X_t)\delta  \tilde X_t\big] +2 \bar K_3(t)+2\l_3H_{pp}(X_t)\delta X_t, \bar \Upsilon_t\big\rangle \\
\dis\qq + \big\langle  -\l_1 \tilde \dbE_{\mathcal{F}_t}\big[ H_{p\mu}(X_t,\tilde X_t)\delta  \tilde X_t \big] +2\l_2 K_4(t)+2\l_3H_{pp}(X_t) \delta X_t,\Upsilon_t\big\rangle \Big].
\ea\right.
\eeaa
Next, by  Assumptions \ref{assum-regH}-(i) and \ref{assum-antidisplacement}-(ii), and \reff{antiH2} we have
\beaa
\left.\ba{lll}
\dis  {d\over dt} \G_t \ge  [4 \underline L^H_{xp} -\l_0 L^H_2]  \dbE[|\bar \Upsilon_t|^2]+\big[2\l_2\underline L^H_{xp}- [|\l_1|+\l_2L_{xx}^V]L^H_2\big] \dbE[|\Upsilon_t|^2]\\
\dis \q +\big[\l_0 \underline L^H_{xx}+2\l_3\underline L^H_{xp}- [ |\l_1|+2\l_3]L^H_2\big]\dbE[|\d X_t|^2]\\
\dis\q  - 2 L^H_2[\l_0 + L^V_{xx}] \dbE[|\Upsilon_t||\bar \Upsilon_t|]\\
\dis\q - \big[|\l_1-2\l_0|L^H_2+2\overline \g  [1+ L^V_{xx}] \underline L^H_{xp} + 2L^V_{xx} L^H_2+2\l_3L_2^H]\big] \left( \dbE[|\d X_t|^2]\right)^{\frac{1}{2}} \left(\dbE[|\bar\Upsilon_t|^2]\right)^{\frac{1}{2}}\\
\dis\q  - L^H_2 \big[|\l_1|+ 2\l_2[1 +L^V_{xx}] +2\l_3\big] \left( \dbE[|\d X_t|^2] \right)^{\frac{1}{2}}\left(\dbE[|\Upsilon_t|^2]\right)^{\frac{1}{2}}.
\ea\right.
\eeaa
Note that, recalling the $\th_1$ in \reff{theta},
\beaa
4 \th_1\dbE[|\bar\Upsilon_t|^2]  +  2\overline \g  [1+ L^V_{xx}]\left( \dbE[|\d X_t|^2]\right)^{\frac{1}{2}}\left( \dbE[|\bar\Upsilon_t|^2] \right)^{\frac{1}{2}} + \th_1 [\l_0\underline \g +2\l_3]\dbE[|\d X_t|^2]  \ge 0,
\eeaa
Then, recalling \reff{A} and denoting $a := \big[ (\dbE[|\bar\Upsilon_t|^2])^{1\over 2},~( \dbE[|\Upsilon_t|^2])^{1\over 2},~   \left(\dbE[|\d X_t|^2] \right)^{\frac{1}{2}}\big]$, 
\beaa
\left.\ba{lll}
\dis  {d\over dt} \G_t \ge  \big[4 [1-\th_1]\underline L^H_{xp} -\l_0 L^H_2]  \dbE[|\bar \Upsilon_t|^2]+\big[2\l_2\underline L^H_{xp}- [|\l_1|+\l_2 L_{xx}^V]L^H_2\big] \dbE[|\Upsilon_t|^2]\\
\dis \q +\big[[1-\th_1][\l_0\underline \g +2\l_3]\underline L^H_{xp}- [ |\l_1|+2\l_3]L^H_2\big]\dbE[|\d X_t|^2]\\
\dis\q  - 2 L^H_2[\l_0 + L^V_{xx}] \dbE[|\Upsilon_t||\bar \Upsilon_t|]\\
\dis\q - L^H_2 \big[|\l_1-2\l_0|+ 2L^V_{xx} +2\l_3\big] \left( \dbE[|\d X_t|^2]\right)^{\frac{1}{2}} \left(\dbE[|\bar\Upsilon_t|^2]\right)^{\frac{1}{2}}\\
\dis\q  - L^H_2 \big[|\l_1|+ 2\l_2[1 +L^V_{xx}] +2\l_3\big] \left( \dbE[|\d X_t|^2] \right)^{\frac{1}{2}}\left(\dbE[|\Upsilon_t|^2]\right)^{\frac{1}{2}}.\\
\dis\q  = a \big[ A_1  \underline L^H_{xp} - A_2 L^H_2\big] a^\top \ge 0,
\ea\right.
\eeaa
where the last inequality thanks to  \reff{Lxp} and the fact that $A_1 \ge 0$.
Thus
\beaa
 (AntiMon)^{\vec \l}_\xi V(0, \eta, \eta)  = \G_0 \le \G_T= (AntiMon)^{\vec \l}_{X_T}G(\d X_T, \d X_T)\le 0.
\eeaa
That is, $V(0,\cd,\cd)$ is $\vec\l$-anti-monotone.
\qed

\section{The Lipschitz continuity}
\label{sect-Lipschitz}
\setcounter{equation}{0}

We first show that the anti-monotonicity of $V$ implies the uniformly Lipschitz continuity of $\pa_x V$ in $\mu$ under $W_2$. Unlike in \cite{GMMZ}, since we do not require the first order derivatives of $G, H$ to be bounded, here we do not expect the Lipschitz continuity of $V$ itself.

\begin{thm}
\label{thm-Lipschitz2}
Let Assumptions \ref{assum-regH}-(i), \ref{assum-regG}-(i) hold and $V$ be a classical solution of the master equation \eqref{master} such that 
$$\pa_{xx}V(t,\cdot,\cdot)\in \cC^2(\mathbb R^d\times\mathcal{P}_2), \quad  \pa_{x\mu}V(t,\cdot,\cdot,\cdot)\in \cC^2(\mathbb R^d\times\mathcal{P}_2\times\mathbb R^d),$$ 
and all the second and higher order derivatives of $V$ involved above are also continuous in the time variable and are uniformly bounded. Assume further that $V(t,\cdot,\cdot)$ is $\vec \l$-anti-monotone  in the sense of \reff{anti} for all $t\in [0, T]$. Then $\pa_x V$ is uniformly Lipschitz continuous in $\mu$ under $W_2$, with  a Lipschitz constant $C^\mu_2$ depending only on $\vec \l$, the parameters in \reff{Hbound} and \reff{Gbound}, and $L^V_{xx}$.
\end{thm}
\proof
In this proof, $C>0$ denotes a generic constant depending only on quantities mentioned in the statement of the theorem. As in the proof of Theorem \ref{thm-antidisplacement}, without loss of generality we show the theorem only for $t_0=0$. First, by \eqref{anti} we have, for any $\xi, \eta\in \dbL^{2}(\cF_t^1)$,
\begin{equation}\label{paxmuVeta}
\dbE\Big[\Big|\mathbb {\tilde E}_{\mathcal{F}_T^1}\big[\pa_{x\mu}V(t,\xi,\mathcal{L}_{\xi},\tilde \xi)\tilde\eta\big]\Big|^2\Big]\leq C\Big|\mathbb {\tilde E}\Big[\big\langle\pa_{x\mu}V(t,\xi,\mathcal{L}_{\xi},\tilde \xi)\tilde\eta,\eta\big\rangle\Big]\Big| +C\mathbb E[|\eta|^2].
\end{equation}
Next, applying H\"older's inequality to \eqref{paxmuVeta} we have
\begin{equation}\label{weakW2Lip}
\dbE\Big[\Big|\mathbb {\tilde E}_{\mathcal{F}_T^1}\big[\pa_{x\mu}V(t,\xi,\mathcal{L}_{\xi},\tilde \xi)\tilde\eta\big]\Big|^2\Big]\leq C\mathbb E[|\eta|^2].
\end{equation}
From now on we fix $\xi\in \dbL^{2}(\cF_0)$ and $\eta\in \dbL^{2}(\cF_0)$ and continue to use the notation as in the proof of Theorem \ref{thm-antidisplacement}. In particular,  $X,\delta X,  \mu_t, \Upsilon, \bar \Upsilon$ are defined by  \reff{XY}. 
Applying \eqref{weakW2Lip} by replacing $\dbE$ with $\dbE_{\cF^0_t}$ and noting that $X_t$ is $\cF_t$-measurable, we have
\bea
\label{deltaX}
\dbE[|\Upsilon_t|^2] =\dbE\Big[\dbE_{\mathcal{F}_t^0}\Big[\Big|\mathbb {\tilde E}_{\mathcal{F}_T}\big[\pa_{x\mu}V(t,X_t,\mu_t,\tilde X_t)\delta\tilde X_t\big]\Big|^2\Big]\Big]\leq C\dbE\left[\dbE_{\mathcal{F}_t^0}[|\delta X_t|^2]\right]\leq C \mathbb E[|\delta X_t|^2].
\eea
Using H\"older's inequality on \eqref{XY} and noting in particular $|\bar \Upsilon_t|\le L^V_{xx} |\d X_t|$, we obtain 
\begin{equation}\label{deltaX1}
|\delta X_t|^2\leq 2|\eta|^2+C\int_0^t\Big[|\delta X_s|^2+ \big| \tilde \dbE_{\mathcal{F}_s}[ |\delta  \tilde X_s|]\big|^2+|\Upsilon_s|^2\Big]ds.
\end{equation}
Taking expectation on \eqref{deltaX1} and using \eqref{deltaX}, we derive
\[
\mathbb E\big[|\delta X_t|^2\big] \leq  2\dbE[|\eta|^2]+ C \int_0^t \dbE\big[|\delta X_s|^2\big]ds.
\]
Then it follows from Gr\"onwall's inequality  that 
\bea\label{EY2}  
 \sup_{t \in [0,T]}\mathbb E\big[|\delta X_t|^2\big] \le C \dbE\big[|\eta|^2\big].
 \eea
 
 Next, by \eqref{dUpsilon}, we have
\beaa
\Upsilon_t = \Upsilon_T - \int_t^T \big[ K_3(s) \Upsilon_s + K_4(s)\big] ds -\int_t^T (dB_s)^\top K_1(s) - \b \int_t^T (dB^0_s)^\top K_2(s).
\eeaa
Take conditional expectation $\tilde \dbE_{\cF_t}$, we have
\bea
\label{Upsilon1}
\Upsilon_t = \tilde \dbE_{\cF_t}\Big[\pa_{x\mu} G(X_T, \mu_T, \tilde X_T) \d\tilde X_T\Big] - \int_t^T \tilde \dbE_{\cF_t}\Big[ K_3(s) \Upsilon_s + K_4(s)\Big] ds.
\eea
Then by \reff{Upsilon1} and the required regularity of $G, H$ and $V$, we have
\beaa
|\Upsilon_t|^2 \le C \tilde \dbE_{\cF_t}\big[|\d\tilde X_T|^2\big] + C\int_t^T \tilde \dbE_{\cF_t}\big[ |\Upsilon_s|^2 + |\d \tilde X_s|^2\big] ds.
\eeaa
Now take conditional expectation $\tilde \dbE_{\cF_0}$, we get
\beaa
\tilde \dbE_{\cF_0}\big[|\Upsilon_t|^2\big] \le C \tilde \dbE_{\cF_0}\big[|\d\tilde X_T|^2\big] + C\int_t^T \tilde \dbE_{\cF_0}\big[ |\Upsilon_s|^2 + |\d \tilde X_s|^2\big] ds.
\eeaa
Thus, by the Gr\"onwall inequality we have
\bea
\label{Upsilon0}
|\Upsilon_0|^2= \tilde \dbE_{\cF_0}\big[|\Upsilon_0|^2\big] \le C \tilde \dbE_{\cF_0}\big[|\d\tilde X_T|^2\big] + C\int_0^T \tilde \dbE_{\cF_0}\big[|\d \tilde X_s|^2\big] ds.
\eea

Note that, recalling the setting in Section \ref{sec:setting}, $\d \tilde X_t$ is measurable with respect to $\cF^0_t \vee \tilde \cF^1_t$, which is independent of $\cF_0$ under $\tilde \dbP$. Then the conditional expectation in the right side of \reff{Upsilon0} is actually an expectation. Plug \reff{EY2} into \reff{Upsilon0}, we have
 \bea
 \label{EUpsilon0}
\Big| \tilde \dbE_{\cF_0}\Big[\pa_{x\mu} V(0, \xi, \mu_0, \tilde \xi) \tilde \eta\Big] \Big|^2 = |\Upsilon_0|^2  \le  C\dbE\big[|\eta|^2\big].
\eea
This implies
\bea
\label{pamuVbound}
\Big|\tilde\dbE\big[\pa_{x\mu} V(0, x,  \mu_0, \tilde \xi)  \tilde \eta \big]\Big|\le  C(\mathbb E|\eta|^2)^{\frac{1}{2}},\q\mu_0-\mbox{a.e.} ~x.
\eea
Since $\pa_\mu V$ is continuous,  then \reff{pamuVbound} actually holds for all $x$. 
In particular, this implies that there exists a constant $C^{\mu_0}_2>0$ such that
\bea\label{eq:w2lip}
\Big|\pa_xV(0, x, \cL_{\xi+\eta}) - \pa_x V(0, x, \cL_\xi)\Big| = \Big|\int_0^1 \dbE\big[\pa_{x\mu} V(0, x, \cL_{\xi+\th \eta}, \xi + \th \eta) \eta\big] d\th\Big|\le C^{\mu_0}_2(\mathbb E[|\eta|^2])^{\frac{1}{2}}.
\eea
Now, taking random variables $\xi,\eta$ such that $W^2_2(\cL_{\xi+\eta},\cL_{\xi})=\mathbb E[|\eta|^2]$, the above inequality
exactly means that $\pa_x V(0,x,\cdot)$ is uniformly Lipschitz continuous in $\mu_0$ under $W_2$ with uniform Lipschitz constant $C^{\mu_0}_2$.
\qed

We emphasize that the above Lipschitz continuity is under $W_2$, while the global wellposedness of the master equation requires the $W_1$-Lipschitz continuity. As in \cite{GMMZ}, we shall derive the desired $W_1$-Lipschitz continuity from the $W_2$-Lipschitz continuity by utilizing the pointwise representation for the Wasserstein derivative developed in \cite{MZ2}. Note again that in Theorem \ref{thm-Lipschitz2} we only have the Lipschitz continuity for $\pa_x V$, but not for $V$, so at below we shall also consider  $\vec U(t,x, \mu) := \pa_x V(t,x,\mu)$, which formally should satisfy the following vectorial master equation on $[0,T)\times\R^d\times\cP_2(\R^d)$, with terminal condition $\vec{U}(T,x,\mu) = \pa_xG(x,\mu)$:
\bea
\label{vecmaster}
\left.\ba{c}
\dis -\pa_t \vec{U} -\frac{\h\b^2}{2} \tr(\pa_{xx} \vec{U}) + \pa_x H(x,\mu,\vec{U})+\pa_pH(x,\mu,\vec{U})\cdot\pa_x \vec{U}  - \vec\cN \vec{U} =0,\q\mbox{where} \\
\dis\vec \cN \vec{U}(t,x,\mu):= \tr\Big(\bar{\tilde  \dbE}\Big[\frac{\h\b^2}{2} \pa_{\tilde x} \pa_\mu \vec{U}(t,x, \mu, \tilde \xi) - \pa_\mu \vec{U}(t, x, \mu, \tilde \xi)(\pa_pH)^\top(\tilde \xi,\mu, \vec{U}(t, \tilde \xi, \mu))
 \\
 +\b^2 \pa_x\pa_\mu \vec{U}(t,x,\mu,\tilde \xi)+\frac{\b^2}{2}\pa_{\mu\mu}\vec{U}(t,x,\mu,\bar\xi,\tilde\xi)\Big]\Big).
 \ea\right.
\eea

To be precise, fix $t_0, \xi$, we first consider the following McKean-Vlasov SDE on $[t_0, T]$: 
\bea
\label{vecFBSDE} 
\left.\ba{lll}
\dis\left.\ba{ll}
\dis X^{\xi}_t & \dis =\xi - \int_{t_0}^t\pa_pH(X_s^\xi,\rho_s,\nabla Y_s^{\xi})ds+B^{t_0}_t+\b B_t^{0,t_0}, \q  \rho_t := \rho^\xi_t:= \cL_{X_t^\xi|\cF^0_t}; \\
\dis \nabla Y_t^{\xi}&\dis =\pa_x G(X_T^{\xi},\rho_T)-\int_t^T\pa_xH(X_s^{\xi},\rho_s,\nabla Y_s^{\xi})ds- \int_t^T\nabla Z_s^{\xi}\cd dB_s-\int_t^T\nabla Z_s^{0,\xi}\cd dB_s^{0}.
\ea\right.
\ea\right.
\eea
Next, given $\rho$ as above, for fixed $x\in \dbR^d$ and letting $(e_1,\cds, e_d)$ denote the natural basis of $\dbR^d$, we introduce a series of FBSDEs, possibly McKean-Vlasov type:
{\small\bea
\label{vecFBSDEx} 
\left\{\ba{lll}
\dis X^{\xi,x}_t  =x - \int_{t_0}^t\pa_pH(X_s^{\xi,x},\rho_s,\nabla Y_s^{\xi,x})ds+B^{t_0}_t+\b B_t^{0,t_0}; \\
\dis \nabla Y_t^{\xi,x} =\pa_x G(X_T^{\xi,x},\rho_T)- \!\int_t^T\! \pa_xH(X_s^{\xi,x},\rho_s,\nabla Y_s^{\xi,x})ds - \int_t^T\nabla Z_s^{\xi,x}\cd dB_s-\int_t^T\nabla Z_s^{0,\xi,x}\cd dB_s^{0};
\ea\right.
\eea
\bea
\label{tdYx}
\left\{
\ba{ll}
\dis \td_{k} X_t^{\xi,x}=\dis e_k -\int_{t_0}^t \Big[(\nabla_k X_s^{\xi,x})^\top\pa_{xp}H(X_s^{\xi,x},\rho_s,\nabla Y_s^{\xi,x})
+ (\nabla_k^2 Y_s^{\xi,x})^\top\pa_{pp}H(X_s^{\xi,x},\rho_s,\nabla Y_s^{\xi,x})\Big]ds; \\
\dis \td_{k}^2 Y_t^{\xi,x}= \dis (\td_{k} X^{\xi,x}_T)^{\top} \pa_{xx} G(X^{\xi,x}_T,\rho_T) -\dis \int_t^T\td_{k}^2 Z_s^{\xi,x}\cdot dB_s^{t_0}- \int_t^T \td_{k}^2 Z_s^{0,\xi,x}\cdot dB_s^{0,t_0} \\
\dis\qq - \dis \int_t^T\Big[ (\td_k X_s^{\xi,x})^{\top} \pa_{xx} H(X_s^{\xi,x},\rho_s,\nabla Y_s^{\xi,x})+(\nabla_k^2 Y_s^{\xi,x})^{\top} \pa_{px}H(X_s^{\xi,x},\rho_s,\nabla Y_s^{\xi,x})\Big]ds; 
\ea\right.
\eea
\bea
\label{tdYx-}
\left\{\ba{ll}
\dis \td_{k} \mathcal X_t^{\xi,x}=\dis-\int_{t_0}^t\Big[(\td_k \mathcal X_s^{\xi,x})^{\top}\pa_{xp}H(X_s^{\xi},\rho_s,\nabla Y_s^\xi)+(\td_k^2 \mathcal Y_s^{\xi,x})^{\top}\pa_{pp}H(X_s^{\xi},\rho_s,\nabla Y_s^\xi)\ss\\
\dis\qq  +\tilde{\mathbb E}_{\mathcal{F}_s}\big[(\nabla_{k} \tilde X_s^{\xi,x})^{\top}(\pa_{\mu p} H)(X^{\xi}_s,\rho_s,\tilde X_s^{\xi,x}, \nabla Y^{\xi}_s)
 +(\nabla_{k} \tilde {\mathcal X}_s^{\xi,x})^{\top}\pa_{\mu p} H(X^{\xi}_s,\rho_s,\tilde X_s^{\xi},\nabla Y^{\xi}_s)\big] \Big]ds; \ms\\
\dis \td_{k}^2 \mathcal Y_t^{\xi,x}= \tilde{\mathbb E}_{\mathcal{F}_T}\Big[(\nabla_{k} \tilde X_T^{\xi,x})^{\top}\pa_{\mu x}G(X_T^\xi,\rho_T,\tilde X_T^{\xi,x})+(\nabla_{k}\tilde {\mathcal X}_T^{\xi,x})^{\top}\pa_{\mu x}G(X_T^\xi,\rho_T, \tilde X_T^\xi)\Big]\ss\\
\dis\qq + ( \td_{k} \mathcal X^{\xi,x}_T)^{\top}\pa_{xx} G(X^{\xi}_T,\rho_T)  -\int_t^T\td_{k}^2 \mathcal Z_s^{\xi,x}\cdot dB_s^{t_0}- \int_t^T \td_{k}^2 \mathcal Z_s^{0,\xi,x}\cdot dB_s^{0,t_0}\\
 \dis\qq - \int_t^T\Big[ (\td_{k} \mathcal X^{\xi,x}_s)^{\top} \pa_{xx} H\big(X_s^{\xi},\rho_s,\nabla Y_s^{\xi})+ ( \td_{k}^2 \mathcal Y^{\xi,x}_s)^{\top} \pa_{px} H\big(X_s^{\xi},\rho_s, \nabla Y_s^{\xi})\\[5pt]
 \dis\qq +\tilde{\mathbb E}_{\mathcal{F}_s}\big[(\nabla_{k} \tilde X_s^{\xi,x})^{\top}\pa_{\mu x}H\big(X^{\xi}_s,\rho_s,\tilde X_s^{\xi,x}, \nabla Y^{\xi}_s )+(\nabla_{k} \tilde {\mathcal X}_s^{\xi,x})^{\top}\pa_{\mu x}H\big(X^{\xi}_s, \rho_s,\tilde X_s^{\xi},\nabla Y^{\xi}_s)\big]\Big]ds;
\ea\right.
\eea
\bea
\label{tdYmu}
\left\{\ba{ll}
\dis \td_{\mu_k} X_t^{x,\xi,\tilde x}=-\int_{t_0}^t\Big[(\td_{\mu_k} X_s^{x,\xi,\tilde x})^{\top}\pa_{xp}H(X_s^{\xi,x},\rho_s,\nabla Y_s^{\xi,x})+(\td_{\mu_k}^2 Y_s^{x,\xi,\tilde x})^{\top}\pa_{pp}H(X_s^{\xi,x},\rho_s,\nabla Y_s^{\xi,x})\ss\\
\dis \q +\tilde{\mathbb E}_{\mathcal{F}_s}\big[(\nabla_{k} \tilde X_s^{\xi,\tilde x})^{\top}\pa_{\mu p} H(X^{\xi,x}_s,\rho_s,\tilde X_s^{\xi,\tilde x}, \nabla Y^{\xi,x}_s)
\dis +(\nabla_{k} \tilde {\mathcal X}_s^{\xi,\tilde x})^{\top}\pa_{\mu p} H(X^{\xi,x}_s,\rho_s,\tilde X_s^{\xi},\nabla Y^{\xi,x}_s)\big] \Big]ds; \ms\\
\td_{\mu_k}^2 Y_t^{x,\xi,\tilde x}=\tilde{\mathbb E}_{\mathcal{F}_T}\Big[(\nabla_{k} \tilde X_T^{\xi,\tilde x})^{\top}\pa_{\mu x} G (X_T^{\xi,x},\rho_T,\tilde X_T^{\xi,\tilde x})+ (\nabla_{k}\tilde {\mathcal X}_T^{\xi,\tilde x})^{\top}\pa_{\mu x} G (X_T^{\xi,x},\rho_T,\tilde X_T^\xi)\Big]\\
\dis\q +(\nabla_{\mu_k}X_T^{x,\xi,\tilde x})^{\top}\pa_{xx}G(X_T^{\xi,x},\rho_T)  -\int_t^T\nabla_{\mu_k}^2 Z_s^{x,\xi,\tilde x}\cdot dB_s- \int_t^T \td_{\mu_k}^2 Z_s^{0,x,\xi,\tilde x}\cdot dB_s^{0}\\
\dis\q - \int_t^T \Big[(\nabla_{\mu_k}X_s^{x,\xi,\tilde x})^{\top}\pa_{xx}H(X_s^{\xi,x},\rho_s,\nabla Y_s^{\xi,x})+(\nabla_{\mu_k}^2 Y_s^{x,\xi,\tilde x})^{\top}\pa_{px}H(X_s^x,\rho_s,\nabla Y_s^{\xi,x})\ss\\
\dis\q +\tilde{\mathbb E}_{\mathcal{F}_s}\big[(\nabla_{k} \tilde X_s^{\xi,\tilde x})^{\top}\pa_{\mu x} H (X_s^{\xi,x},\rho_s,\tilde X_s^{\xi,\tilde x},\nabla Y_s^{\xi,x})+(\nabla_{k}\tilde {\mathcal X}_s^{\xi,\tilde x})^{\top}\pa_{\mu x} H (X_s^{\xi,x},\rho_s,\tilde X_s^\xi,\nabla Y_s^{\xi,x})\big]\Big]ds.
\ea\right.
\eea}

The following local (in time) result provides the crucial $W_1$-Lipschitz continuity of $\vec{U}$.

\begin{prop}\label{prop:pamuV}  Let Assumptions \ref{assum-regH}-(i) and \ref{assum-regG}-(i) hold. Recall the constants $\overline L_{xx}^{H},\overline L_{xp}^H,L_2^{H}$ in \reff{Hbound},  $L^G_2$, $\overline L^G_{xx}$ in \reff{Gbound}, and $\tilde L_2^G$ in Remark \ref{rem:W2Lip}.  Then there exists a constant $\d>0$, depending only $d$, $\overline L_{xx}^{H}$, $\overline L_{xp}^H$, $L_2^{H}$, $\overline L_{xx}^G$, $\tilde L_2^G$,  such that whenever $T-t_0\le \d$, 
the following hold.

(i) The McKean-Vlasov FBSDEs \eqref{vecFBSDE}, \eqref{vecFBSDEx}, \eqref{tdYx}, \eqref{tdYx-}, and \reff{tdYmu} are well-posed on $[t_0,T]$, for any $\mu \in \mathcal P_2$ and $\xi \in \mathbb L^2(\mathcal F_{t_0}, \mu)$. 

(ii) Define $\vec{U}(t_0,x, \mu):= \nabla Y_{t_0}^{x,\xi}$. Then we have the pointwise representation:
\bea\label{pamuV}
\pa_{\mu_k} \vec{U}(t_0,x,\mu,\tilde x)=\nabla_{\mu_k}^2Y_{t_0}^{x,\xi,\tilde x}.
\eea
Moreover, there exists a constant $C_1^\mu >0$, depending only on $d, \overline L^{H}_{xx}, \overline L_{xp}^H,  L^{H}_2, L_2^G,\overline L_{xx}^G$ such that
\begin{equation}\label{W1Lip}
|\pa_{\mu}\vec{U}(0,x,\mu,\tilde x)|\leq C_1^\mu.
\end{equation}

(iii) Assume further that Assumptions  \ref{assum-regH}-(ii) and \ref{assum-regG}-(ii) hold true.  Then the vectorial master equation \reff{vecmaster} has a unique classical solution $\vec U$. Moreover, 
$$
\vec{U}(t,\cdot,\cdot),\;\; \pa_{x}\vec{U}(t,\cdot,\cdot)\in \cC^2(\mathbb R^d\times\mathcal{P}_2),\;\; \pa_{\mu}\vec{U}(t,\cdot,\cdot,\cdot)\in \cC^2(\mathbb R^d\times\mathcal{P}_2\times\mathbb R^d),
$$ 
and all their derivatives in the state and probability measure variables are continuous in the time variable and are uniformly bounded.

(iv) The following decoupled McKean-Vlasov FBSDE 
\bea
\label{masterFBSDE} 
\left.\ba{lll}
\dis X^{x}_t =x+B^{t_0}_t+\b B_t^{0,t_0}; \\
\dis Y_t^{x,\xi}=G(X_T^{x},\rho_T)-\int_t^TH(X_s^{x},\rho_s,\vec{U}(s,X_s^{x},\rho_s))ds- \int_t^T Z_s^{x,\xi}\cd dB_s-\int_t^TZ_s^{0,x,\xi}\cd dB_s^{0}
\ea\right.
\eea
is well-posed on $[t_0,T]$ for any $x\in\mathbb R^d$. Define $V(t_0,x,\mu):=Y_{t_0}^{x,\xi}$. Then $V$ is the unique classical solution of the master equation \eqref{master}  and $\pa_xV=\vec{U}$ on $[0,T]\times\mathbb R^d\times\cP_2$.
\end{prop}
We emphasize that at above $C^\mu_1$ depends on $L^G_2$ in \reff{Gbound}, but the $\d$ depends only on $\tilde L^G_2$ in  Remark \ref{rem:W2Lip}, not on $L^G_2$. 

\ms
\begin{proof}
The proof of (i)-(iii) is very lengthy, but essentially identical to as  that of  \cite[Proposition 6.2]{GMMZ}, except that \cite{GMMZ} considers both $\pa_\mu V$ and $\pa_{x \mu} V = \pa_\mu \vec U$. Therefore we omit it here. 

(iv) By the smoothness of  $\vec U$ obtained in (iii), clearly the $V$ defined in (iv) is smooth and $Y^{x,\xi}_t= V(t, X^x_t, \rho_t)$.  By applying It\^{o}'s formula \reff{Ito} we see that $V$ satisfies the PDE:
\bea
\label{VUPDE}
\left.\ba{c}
\dis-\pa_t V -\frac{\h\b^2}{2} \tr(\pa_{xx} V) + H(x,\mu,\vec U)  - \tr\bigg(\bar{\tilde  \dbE}\Big[\frac{\h\b^2}{2} \pa_{\tilde x} \pa_\mu V(t,x, \mu, \tilde \xi) 
+\frac{\b^2}{2}\pa_{\mu\mu}V(t,x,\mu,\bar\xi,\tilde\xi) \\
 \dis - \pa_\mu V(t, x, \mu, \tilde \xi)(\pa_pH)^\top(\tilde \xi,\mu,  \vec U(t, \tilde \xi, \mu))+\b^2 \pa_x\pa_\mu V(t,x,\mu,\tilde \xi)\Big]\bigg)=0.
 \ea\right.
\eea
Differentiate it with respect to $x$,  we obtain the PDE for $\vec U' := \pa_x V$:
\bea
\label{VUPDE2}
\left.\ba{lll}
\dis -\pa_t \vec{U'} -\frac{\h\b^2}{2} \tr(\pa_{xx} \vec{U'}) + \pa_x H(x,\mu,\vec{U})+\pa_pH(x,\mu,\vec{U})\cdot\pa_x \vec{U}  \\
\dis- \tr\Big(\bar{\tilde  \dbE}\Big[\frac{\h\b^2}{2} \pa_{\tilde x} \pa_\mu \vec{U'}(t,x, \mu, \tilde \xi) +\frac{\b^2}{2}\pa_{\mu\mu}\vec{U'}(t,x,\mu,\bar\xi,\tilde\xi)\ms\\
\dis\qq - \pa_\mu \vec{U'}(t, x, \mu, \tilde \xi)(\pa_pH)^\top(\tilde \xi,\mu, \vec{U}(t, \tilde \xi, \mu))+\b^2 \pa_x\pa_\mu \vec{U'}(t,x,\mu,\tilde \xi)\Big]\Big) =0.
 \ea\right.
\eea
Compare this with \reff{vecmaster}, we see that $\vec U$ also satisfies \reff{VUPDE2}. Thus, by the uniqueness we have $\vec U = \vec U' = \pa_x V$. Plug this into \reff{VUPDE} we verify that $V$ satisfies \reff{master}.
\end{proof}

\section{Uniform estimates of $\pa_{xx}V$}
\label{sect-Vxx}
\setcounter{equation}{0}

We note that all the above results rely on the bound $L^V_{xx}$ of $\pa_{xx} V$ in \reff{Vxx}. In particular, in Theorem \ref{thm-antidisplacement} the $\underline L^H_{xp}$ depends on $L^V_{xx}$. Then it is crucial to obtain an a priori uniform estimate of $L^V_{xx}$ which is independent of $\underline L_{xp}^{H}$. Recall \reff{YXV}, we have $\pa_{xx} V = \pa_{xx} u$, so it suffices to establish the a priori estimate for the solution $u$ to the backward SPDE in \reff{SPDE}, for an arbitrarily given $\rho$ (not necessarily satisfying the forward SPDE in \reff{SPDE}).
 
For this purpose we consider a special form of $H$. 
\begin{assum} \label{assum-regH1}  $H$ takes the following form: 
\begin{equation}\label{eq:Hform}
H(x,\mu,p)=\langle A_0 x, p\rangle+H_0(x,\mu,p)
\end{equation}
where $A_0\in\mathbb R^{d\times d}$ is a constant matrix and $H_0: \mathbb R^d\times\mathcal{P}_2\times\mathbb R^d\to\mathbb R$ is a function satisfying

(i) $H_0\in \cC^2(\mathbb R^d\times\mathcal{P}_2\times\mathbb R^d)$ and there exist constants $\underline L_{xx}^{H_0}, \overline L_{xx}^{H_0}$,$L_2^{H_0}>0$ such that 
\bea
\label{H0bound1}
&\dis  \underline \k(\pa_{xx}H_0) \ge \underline L_{xx}^{H_0},\q |\pa_{xx}H_0|\leq \overline L^{H_0}_{xx},\\
\label{H0bound2}
&\dis |\pa_{xp}H_0|,|\pa_{pp}H_0|,|\pa_{x\mu}H_0|,|\pa_{p\mu}H_0|\le L_2^{H_0}~ {\rm{on}}\ \mathbb R^d\times\mathcal{P}_2\times\mathbb R^d.
\eea

(ii) $H_0$ satisfies Assumption \ref{assum-regH}-(ii). 
\end{assum}

Given $A_0$, consider its Jordan decomposition:
\begin{equation}
\label{Jordan}
A_0=Q_0  J_0 Q_0^{-1}, 
\end{equation}
where $J_0\in \mathbb C^{d\times d}$ is the Jordan normal form of $A_0$ and $Q_0\in \dbC^{d\times d}$ is invertible. Let $\bar Q_0$ denote the conjugate of $Q_0$ and thus $Q_0 \bar Q_0^\top$ is  positive definite. The following estimate will be crucial.  

\begin{lem}
\label{lem-Aest} Recall \reff{kappaA}. For any $t\ge 0$, we have
\bea
\label{Aest}
|e^{-A_0 t}| \le \sqrt{L^{A_0}} e^{[1-\underline\k'(A_0)]t},\q \mbox{where}\q L^{A_0} := \inf_{Q_0} ~ {\overline \k(Q_0 \bar Q_0^\top)\over \underline \k(Q_0 \bar Q_0^\top)}.
\eea
Here the infimum is over all $Q_0$ satisfying \reff{Jordan}. 
\end{lem}
\proof Fix $J_0, Q_0$ as in \reff{Jordan}. It is obvious that $e^{-A_0 t} = Q_0 e^{-J_0 t} Q_0^{-1}$. We claim that
\bea
\label{J0claim}
\big|\langle e^{-J_0 t}x, y\rangle\big|\le e^{[1-\underline\k'(A_0)]t} |x||y|,\q\forall x, y\in \dbC^d.
\eea
Then, for any $x, y\in \dbR^d$ with $|x|=|y|=1$, we have
\beaa
&&\dis \Big|\big\langle e^{-A_0 t}x, y\big\rangle\Big| = \Big|\big\langle e^{-J_0 t} Q_0^{-1}x, Q_0\top y\big\rangle\Big|\le e^{[1-\underline\k'(A_0)]t} |Q_0^{-1}x||Q_0\top y|\\
&&\dis \le e^{[1-\underline\k'(A_0)]t} \sqrt{\overline \k(Q_0^{-1}(\bar Q_0^\top)^{-1})}\sqrt{ \overline \k(Q_0\bar Q_0^\top)} = e^{[1-\underline\k'(A_0)]t} \sqrt{ \overline \k(Q_0 \bar Q_0^\top)\over \underline \k(Q_0 \bar Q_0^\top)}.
\eeaa
Since $Q_0$ is arbitrary, this implies \reff{Aest} immediately. 

To see \reff{J0claim}, assume the Jordan normal form $J_0 = diag(J_1, \cds, J_k)$. Here $d_1 +\cds + d_k=d$; $J_i = \l_i I_{d_i} + U_{d_i}\in \dbR^{d_i\times d_i}$, $i=1,\cds, k$; $\l_1,\cds, \l_k$ are all the eigenvalues of $A_0$; and $U_{d_i}$ is the matrix whose $(j, j+1)$-component is $1$, $j=1,\cds, d_i-1$, and all other components are $0$. It is straightforward to see that
\beaa
e^{-J_0 t} = diag(e^{-J_1 t}, \cds, e^{-J_kt}).
\eeaa
Note that, for each $i$, since $I_{d_i}$ and $U_{d_i}$ can commute, and $U_{d_i}^{d_i}=0$,
\beaa
e^{-J_i t} = e^{-\l_it} e^{-U_{d_i} t } = e^{-\l_i t} \sum_{n=0}^{d_i-1}{(-t)^n\over n!} U_{d_i}^n.
\eeaa
For any $x^{(i)}, y^{(i)} \in \dbC^{d_i}$, it is clear that 
\beaa
\Big|\langle U_{d_i}^n x^{(i)}, y^{(i)}\rangle\Big| \le {1\over 2}[|x^{(i)}|^2+|y^{(i)}|^2].
\eeaa 
Then, for $x = (x^{(1)},\cds, x^{(k)}),  y = (y^{(1)},\cds, y^{(k)})\in \dbC^d$ with $|x|=|y|=1$, we have
\beaa
&&\dis \Big|\langle e^{-J_0 t} x, y\rangle\Big| = \Big|\sum_{i=1}^k \langle e^{-J_i t} x^{(i)}, y^{(i)}\rangle\Big| \le \sum_{i=1}^k |e^{-\l_i t}|  \sum_{n=0}^{d_i-1}{t^n\over n!} \Big|\langle U_{d_i}^n x^{(i)}, y^{(i)}\rangle\Big| \\
&&\dis \le e^{-\underline\k'(A_0)t}\sum_{i=1}^k   \sum_{n=0}^{d_i-1}{t^n\over n!} {1\over 2}[|x^{(i)}|^2+|y^{(i)}|^2] \le e^{-\underline\k'(A_0)t}\sum_{n=0}^{d-1}{t^n\over n!}.
\eeaa
This implies \reff{J0claim} immediately.
\qed

\begin{rem}\label{rem:pxxV}
(i) The form \eqref{eq:Hform} is assumed for the estimate \reff{Aest} and for the property 
\bea
\label{commute}
d e^{-A_0 t} = - e^{-A_0t} A_0 dt = -A_0 e^{-A_0 t} dt,
\eea
 required in the proof of Theorem  \ref{thm:xxbdd} below. In general $e^{-\int_0^t \pa_{xp} H ds}$ does not enjoy these properties. When $d=1$, however, $e^{-\int_0^t \pa_{xp} H ds}$ obviously satisfies similar properties and thus we do not need the special form  \eqref{eq:Hform}. Moreover, we remark that any alternative structures which could ensure a uniform a priori bound for $\pa_{xx} u$ can serve our purpose.

(ii) It is clear that, under \reff{eq:Hform}, \reff{H0bound1}, and \reff{H0bound2},  we may set 
\bea
\label{HAbound}
\underline L^H_{xp} := \underline \k(A_0) - L^{H_0}_2,\q  \overline L^H_{xp} := |A_0|   + L^{H_0}_2;\q  \underline L^H_{xx} := \underline L^{H_0}_{xx},\q \overline L^H_{xx}:=\overline L^{H_0}_{xx},\q L^H_2:= L^{H_0}_2.
\eea
Then \reff{Hbound} and \reff{antiH1} hold true. We shall remark though that the term $\underline \k(A_0)$ and the condition $\underline \k(\pa_{xx}H_0) \ge \underline L_{xx}^{H_0}$ are not used in Theorem \ref{thm:xxbdd} below.

(iii) When $A_0$ is symmetric, one can easily see that $L^{A_0} = 1$, and in this case \reff{Aest}  can be  improved:
$|e^{-A_0 t}| \le e^{-\underline\k'(A_0)t}$.
\end{rem}

Then we have the following uniform a priori estimate. 
 
\begin{thm}
\label{thm:xxbdd}
Let Assumptions \ref{assum-regG}-(i), \ref{assum-regH1} hold and $\rho:[0,T]\times\Omega\to\mathcal{P}_2$ be $\mathbb F^0$-progressively measurable with $\dis\sup_{t\in[0,T]}\mathbb E\big[\int_{\dbR^d}|x|^2\rho_t(dx)\big]<+\infty$.    Assume $(u, v)$ is a classical solution to the backward SPDE in \reff{SPDE} for the given $\rho$ here ($\rho$ is not necessarily a solution to the forward SPDE in \reff{SPDE}) such that $\pa_{xx} u$ is bounded and, for some fixed constant $\overline L^{A}\ge 1$,
\bea
\label{LA}
\left.\ba{c}
L^{A_0} \le \overline L^{A},\q \underline \k'(A_0)  \ge \th_2 :=\max\Big\{\th_3,~{\overline L_{xx}^{H_0}\over 2\overline L_{xx}^G}+1\Big\},\ms\\
\mbox{where}\q \th_3 := 1+L_2^{H_0}\overline L^{A} \Big[1+\overline L_{xx}^G\overline L^{A}+\sqrt{(1+\overline L_{xx}^G\overline L^{A})^2-1}~\Big].
\ea\right.
\eea
 Then the following estimate holds: 
\bea
\label{pauest1}
\left.\ba{c}
\dis |\pa_{xx}  u(t,x)| \le L_{xx}^u(\th_3),\q\forall (t,x),\q  \mbox{where}\ms\\
\dis L_{xx}^u(\th):=\frac{\th -1 - L_2^{H_0}\overline L^{A}-\sqrt{(\th-1 -L_2^{H_0}\overline L^{A})^2-2L_2^{H_0}\overline L_{xx}^G(\overline L^{A})^2 [\th-1] }}{L_2^{H_0}\overline L^{A}}, ~ \th \ge \th_3.
\ea\right.
 \eea
\end{thm}
We note that \reff{LA} implies $L^u_{xx}(\th)$ is well-defined for $\th\ge \th_3$, and we emphasize that the bound $L^u_{xx}(\th_3)$ depends only on $L^{H_0}_2$, $\overline L^G_{xx}$ and $\overline L^{A}$, in particular not on $T$, $\underline \k'(A_0)$, or $\overline L^{H_0}_{xx}$.

\ms
\proof Fix $(t_0, x)$. First, under our conditions it is clear that the following FBSDE on $[t_0,T]$ has a unique solution $(X^{x},\nabla Y^{x},\nabla Z^{x},\nabla Z^{0,x})$:
\bea
\label{vecFBSDEx2} 
\left.\ba{lll}
\dis X^{x}_t =x - \int_{t_0}^t\pa_pH(X_s^{x},\rho_s,\nabla Y_s^{x})ds+B^{t_0}_t+\b B_t^{0,t_0}; \\
\dis \nabla Y_t^{x}= \pa_x G(X_T^{x},\rho_T)-\int_t^T\pa_xH(X_s^{x},\rho_s,\nabla Y_s^{x})ds- \int_t^T\nabla Z_s^{x}\cd dB_s-\int_t^T\nabla Z_s^{0,x}\cd dB_s^{0}.
\ea\right.
\eea
In particular, $\pa_x u$ serves as the decoupling field:
\bea
\label{decoupling}
\td Y^x_t = \pa_x u(t, X^x_t),\q t\in [t_0, T].
\eea

Next, denote $L_0 := L^u_{xx}(\underline \k'(A_0))$, and consider the following BSDE on $[t_0, T]$:
\bea\label{BSDE-x}
\left.\ba{lll}
\dis \nabla ^2Y_t^x=\pa_{xx}G(X_T^x,\rho_T) -\int_t^T\nabla^2 Z_s^x\cdot dB_s-\int_t^T\nabla^2Z_s^{0,x}\cdot dB_s^0\\
\dis\qq - \int_t^T\Big[\nabla^2 Y_s^x\big[A_0^{\top}+\pa_{px}H_0(X_s^x,\rho_s,\nabla Y_s^{x})\big] +\big[A_0+\pa_{xp}H_0(X_s^x,\rho_s,\nabla Y_s^{x})\big]\nabla^2 Y_s^x\ms\\
\dis\qq +\pa_{xx}H_{0}(X_s^x,\rho_s,\nabla Y_s^{x}) +\big[\nabla^2 Y_s^x \wedge L_0\big]\pa_{pp}H_0(X_s^x,\rho_s,\nabla Y_s^{x})\big[\nabla^2 Y_s^x\wedge L_0\big]\Big]ds.
\ea\right.
\eea
Here $A \wedge L_0 := [ (-L_0) \vee a_{ij} \wedge L_0]_{i, j}$ is the truncated matrix.  The above BSDE has a Lipschitz continuous driver and thus is  well-posed. Recalling \reff{commute} and applying It\^{o}'s formula we have
\beaa
\left.\ba{lll}
\dis e^{-A_0 t} \nabla^2 Y_t^xe^{-A_0^\top t}=e^{-A_0 T}\pa_{xx}G(X_T^x,\rho_T) e^{-A_0^\top T}-\int_t^Te^{-A_0 s}\Big[\nabla^2 Z_s^x\cdot dB_s+\nabla^2Z_s^{0,x}\cdot dB_s^0\Big]e^{-A_0 s}\\
\dis\qq - \int_t^Te^{-A_0 s}\Big[\nabla^2 Y_s^x\pa_{px}H_0(X_s^x,\rho_s,\nabla Y_s^{x}) +\pa_{xp}H_0(X_s^x,\rho_s,\nabla Y_s^{x})\nabla^2 Y_s^x\ms\\
\dis\qq +\pa_{xx}H_{0}(X_s^x,\rho_s,\nabla Y_s^{x}) +\big[\nabla^2 Y_s^x \wedge L_0\big]\pa_{pp}H_0(X_s^x,\rho_s,\nabla Y_s^{x})\big[\nabla^2 Y_s^x\wedge L_0\big]\Big]e^{-A_0^\top s}ds.
\ea\right.
\eeaa
Take conditional expectation $\mathbb E_{\mathcal{F}_{t}}$ on both sides,  we obtain
\beaa
&&\dis \nabla^2 Y_t^x=e^{A_0(t-T)}\mathbb E_{\mathcal{F}_{t}}\big[\pa_{xx}G(X_T^x,\rho_T)\big]e^{A^\top_0(t-T)}\\
&&\dis\qq -\int_t^{T}e^{A_0(t-s)}\mathbb E_{\mathcal{F}_{t}}\Big[\nabla^2 Y_s^x\pa_{px}H_0(X_s^x,\rho_s,\nabla Y_s^{x}) +\pa_{xp}H_0(X_s^x,\rho_s,\nabla Y_s^{x})\nabla^2 Y_s^x\ms\\
&&\dis\qq +\pa_{xx}H_{0}(X_s^x,\rho_s,\nabla Y_s^{x}) +\big[\nabla^2 Y_s^x \wedge L_0\big]\pa_{pp}H_0(X_s^x,\rho_s,\nabla Y_s^{x})\big[\nabla^2 Y_s^x\wedge L_0\big]\Big] e^{A^\top_0(t-s)}ds.
\eeaa
Recall \reff{kappaA2} and apply Lemma \ref{lem-Aest}, we have
\beaa
&&\dis \big| \nabla^2 Y_t^x\big|\leq e^{2[1-\underline \k'(A_0)] (T-t)}\overline L_{xx}^G \overline L^{A}+\frac{\overline L_{xx}^{H_0}\overline L^{A}}{2[\underline \k'(A_0)-1] }[1-e^{2[1-{\underline \k'(A_0)}] (T-t)}]\\
&&\dis\qq +L_2^{H_0}\overline L^{A}[2+L_0]\int_t^Te^{2[1-{\underline \k'(A_0)}] (s-t)}\mathbb E_{\mathcal{F}_{t}}[|\nabla^2Y_s^x|]ds.
\eeaa
Taking the conditional expectation $\mathbb E_{\mathcal{F}_{t_0}}$ and noting that ${\underline \k'(A_0)} \ge \th_2\ge {\overline L_{xx}^{H_0}\over 2\overline L_{xx}^G}+1$, we derive
\begin{eqnarray*}
\mathbb E_{\mathcal{F}_{t_0}}\big[\big|\nabla^2Y_{t}^x\big|\big]&\leq& e^{2[1-\underline \k'(A_0)] (T-t)}\overline L_{xx}^G \overline L^{A}+\overline L_{xx}^G \overline L^{A}[1-e^{2[1-{\underline \k'(A_0)}] (T-t)}]\\
&& +L_2^{H_0}\overline L^{A}[2+L_0]\int_t^Te^{2[1-{\underline \k'(A_0)}] (s-t)}\mathbb E_{\mathcal{F}_{t_0}}[|\nabla^2Y_s^x|]ds\\
&\leq&\overline L_{xx}^G\overline L^{A}+L_2^{H_0}\overline L^{A}[2+L_0]\int_t^Te^{2[1-{\underline \k'(A_0)}] (s-t)}\mathbb E_{\mathcal{F}_{t_0}}[|\nabla^2Y_s^x|]ds.
\end{eqnarray*}
Then by Gr\"onwall's inequality we have
\bea
\label{Gronwall}
\left.\ba{c}
\dis\mathbb E_{\mathcal{F}_{t_0}}\big[\big|\nabla^2Y_{t}^x\big|\big] \le \overline L_{xx}^G\overline L^{A} + \frac{\overline L_{xx}^GL_2^{H_0}|\overline L^{A}|^2[2+L_0]}{{2[\underline \k'(A_0)}-1] -L_2^{H_0}\overline L^{A}[2+L_0]}\times\ms\\
\dis\Big[1- e^{-\big[2[\underline \k'(A_0)-1] -L_2^{H_0}\overline L^{A}[2+L_0]\big] [T-t]}\Big].
\ea\right.
\eea
Recall  \reff{pauest1}, one can check straightforwardly that 
\bea
\label{Luxxdecreasing}
{d\over d\th} L_{xx}^u(\th)={1\over  L^{H_0}_2\overline L^{A}}\Big[1 - {(\th-1 -L_2^{H_0}\overline L^{A})-L_2^{H_0}\overline L_{xx}^G(\overline L^{A})^2 \over \sqrt{(\th-1 -L_2^{H_0}\overline L^{A})^2-2L_2^{H_0}\overline L_{xx}^G(\overline L^{A})^2 [\th-1] }}\Big] <0, ~\forall \th \ge \th_3.
 \eea
Then, since $\underline \k'(A_0) \ge \th_2\ge \th_3$ and $L_0=L^u_{xx}(\underline\k'(A_0))$,  by \reff{LA} and \reff{pauest1} we have
\beaa
2[\underline \k'(A_0)-1] -L_2^{H_0}\overline L^{A}[2+L_0]\ge 2[\th_3-1] -L_2^{H_0}\overline L^{A}[2+L^u_{xx}(\th_3)] \ge 0.
\eeaa 
Thus \reff{Gronwall} implies
\beaa
\dis \mathbb E_{\mathcal{F}_{t_0}}\big[\big|\nabla^2Y_{t}^x\big|\big] &\le& \overline L_{xx}^G\overline L^{A} + \frac{\overline L_{xx}^GL_2^{H_0}|\overline L^{A}|^2[2+L_0]}{{2[\underline \k'(A_0)}-1] -L_2^{H_0}\overline L^{A}[2+L_0]} \\
\dis &=& \frac{2\overline L_{xx}^G\overline L^{A}[\underline \k'(A_0)-1]}{{2[\underline \k'(A_0)-1]} -L_2^{H_0}\overline L^{A}[2+L_0]} = L_0,
\eeaa
where the last equality is due to the straightforward calculation. In particular, by setting $t=t_0$, we have $\big|\nabla^2Y_{t_0}^x\big| \le L_0$. Similarly we can show $\big|\nabla^2Y_{t}^x\big| \le L_0$ for all $t\in [t_0, T]$.  Then  $\nabla^2 Y_s^x\wedge L_0 = \nabla^2 Y_s^x$ and thus \reff{BSDE-x} becomes
\bea\label{BSDE-x2}
\left.\ba{lll}
\dis \nabla ^2Y_t^x=\pa_{xx}G(X_T^x,\rho_T) -\int_t^T\nabla^2 Z_s^x\cdot dB_s-\int_t^T\nabla^2Z_s^{0,x}\cdot dB_s^0\\
\dis\qq - \int_t^T\Big[\nabla^2 Y_s^x\big[A_0^{\top}+\pa_{px}H_0(X_s^x,\rho_s,\nabla Y_s^{x})\big] +\big[A_0+\pa_{xp}H_0(X_s^x,\rho_s,\nabla Y_s^{x})\big]\nabla^2 Y_s^x\ms\\
\dis\qq +\pa_{xx}H_{0}(X_s^x,\rho_s,\nabla Y_s^{x}) +\nabla^2 Y_s^x \pa_{pp}H_0(X_s^x,\rho_s,\nabla Y_s^{x})\nabla^2 Y_s^x\Big]ds.
\ea\right.
\eea

By considering the equation for $\pa_{xx} u$ derived from the BSPDE in \reff{SPDE}, one can readily see from \reff{vecFBSDEx2}, \reff{decoupling}, and \reff{BSDE-x2} that $\nabla ^2Y_t^x = \pa_{xx} u(t, X^x_t)$. In particular, $|\pa_{xx} u(t_0, x)| = |\nabla ^2Y_{t_0}^x|\le L_0$. Since $(t_0, x)$ is arbitrary, we have $|\pa_{xx} u(t, x)| \le L_0 = L^u_{xx}(\underline \k'(A_0))$ for all $(t,x)$.  This, together with \reff{Luxxdecreasing}, implies \reff{pauest1}.
\qed

\section{Global well-posedness}
\label{sect-wellposedness}
\setcounter{equation}{0}

In this section we establish the global well-posedness of the master equation.  We shall first construct the global well-posedness of the vectorial master equation \eqref{vecmaster}. Following the idea in \cite{CCD, CD2,MZ2,GMMZ}, the key is to extend a local classical solution to a global one through an a priori uniform Lipschitz continuity estimate of the solution in $\mu$. We note that Theorem \ref{thm:xxbdd} implies the uniform a priori bound of $\pa_{xx} V$. Then, by applying Theorem \ref{thm-antidisplacement} and \ref{thm-Lipschitz2}, we obtain the uniform a priori Lipschitz continuity of $\vec{U}=\pa_{x}V$ with respect to $\mu$ under $W_2$. Moreover, by Proposition \ref{prop:pamuV} we derive the desired uniform a priori Lipschitz continuity of $\vec{U}$ with respect to $\mu$ under $W_1$.

We now present the main well-posedness result. Note that the dependence on the parameters is quite subtle, so we will introduce them carefully. Following the order of the assumptions below, one can easily construct a class of $G$ and $H$ satisfying all of them, see e.g.  Example \ref{eg:allassum}. 
In particular, in light of Lemma \ref{rem:pxxV} (iii), we may set $\overline L^A = 1$ and consider symmetric $A_0$.

\begin{thm}\label{thm:wellposedness}
Let Assumption \ref{assum-regG} with $\overline L^G_{xx}, L^G_2$ and Assumption \ref{assum-antidisplacement} (i) with $\vec \l\in D_4$ hold true, and $H$ takes the form \reff{eq:Hform} such that Assumption \ref{assum-regH1} (ii) holds and there exists $L^{H_0}_2$ satisfying the requirements in \reff{H0bound2}. Fix an arbitrary $\overline L^A \ge 1$ and set $\th_3$ as in \reff{LA} and $L^V_{xx} := L^u_{xx}(\th_3)$  as in \reff{pauest1}. Assume further the following hold true.

(i) There exist $0<\underline \g<\overline \g$ such that $\underline \g \le \overline L^G_{xx}$, $\overline \g >1$, and \reff{theta} holds true. 

(ii) Set $A_1, A_2$ as in \reff{A}. The matrix $A_0$ satisfies:
\bea
\label{kA0}
L^{A_0} \le \overline L^A,\q \underline\k(A_0) \ge [1+\underline \k(A_1^{-1}A_2)] L^{H_0}_2,\q \underline\k'(A_0)\geq \th_3 ,\q |A_0| + L^{H_0}_2 \le \overline \g[\underline \k(A_0)-L^{H_0}_2].
\eea

(iii) There exist $0<\underline L^{H_0}_{xx}\le \overline L^{H_0}_{xx}$ satisfying \reff{H0bound1} and
\bea
\label{LH0}
\underline \g[\underline \k(A_0)-L^{H_0}_2] \le \underline L^{H_0}_{xx} \le \overline L^{H_0}_{xx} \le \Big[\overline \g[\underline \k(A_0)-L^{H_0}_2]\Big] \wedge \Big[2\overline L^G_{xx}[\underline \k'(A_0)-1]\Big].  
\eea

\no Then the master equation \eqref{master} on $[0, T]$ admits a unique classical solution $V$ with bounded $ \pa_xV$, $\pa_{xx} V$ and $\pa_{x\mu} V$. 

Furthermore, the McKean-Vlasov FBSDEs \eqref{vecFBSDE}, \eqref{vecFBSDEx}, \eqref{tdYx}, \eqref{tdYx-}, \reff{tdYmu} and \reff{masterFBSDE} are also well-posed on $[0,T]$  and the representation formula \reff{pamuV} remains true on $[0, T]$. 
\end{thm}
\proof The uniqueness as well as the wellposedness of the involved FBSDEs and the representation formula \reff{pamuV}  follow exactly the same arguments as in \cite[Theorem 6.3]{GMMZ}. Thus we shall only prove the existence.

Set $\underline L^H_{xp}, \overline L^H_{xp},  \underline L^H_{xx}, \overline L^H_{xx}, L^H_2$ as in \reff{HAbound}. Then clearly Assumptions \ref{assum-regH} and \ref{assum-antidisplacement} hold true. By \reff{kA0} and \reff{LH0} we see that \reff{LA} holds true and thus we have the a priori estimate \reff{pauest1}. Moreover,  by \reff{kA0} we have $\underline L^H_{xp} \ge \underline \k(A_1^{-1}A_2) L^H_2$, and thus the result of Theorem \ref{thm-antidisplacement} holds true.

 We now let $C^\mu_2$ be the a priori (global) uniform Lipschitz estimate of $\pa_xV$ with respect to $\mu$ under $W_2$, as established by Theorems \ref{thm-antidisplacement} and \ref{thm-Lipschitz2}. Let $\d>0$ be the constant in Proposition \ref{prop:pamuV}, but with $\overline L^G_{xx}$ replaced with $L_{xx}^V$ and $L^G_2$ replaced with $C^\mu_2$. Let $0=T_0<\cds<T_n=T$ be a partition such that $T_{i+1}-T_i \le {\d\over 2}$, $i=0,\cds, n-1$. 
 
First, since $T_n - T_{n-2}\le \d$, by Proposition \ref{prop:pamuV}  the master equation \reff{master}  on $[T_{n-2}, T_n]$ with terminal condition $G$ has a unique classical solution $V$. For each $t\in [T_{n-2}, T_n]$, applying Theorem \ref{thm:xxbdd} we have $|\pa_{xx}V(T_{n-1},\cdot,\cdot)|\le L_{xx}^V$. Note that by Proposition \ref{prop:pamuV}-(iii)(iv) $V(t,\cdot,\cdot)$ has further regularities, this enables us to apply Theorems \ref{thm-antidisplacement} and \ref{thm-Lipschitz2} and obtain that $\pa_xV(t,\cdot,\cdot)$ is uniform Lipschitz continuous in $\mu$ under $W_2$ with Lipschitz constant $C^\mu_2$. Moreover, by Proposition \ref{prop:pamuV}-(ii) $\pa_x V(T_{n-1},\cdot,\cdot)$ is also uniformly Lipschitz continuous in $\mu$ under $W_1$.

We next consider  the master equation \reff{master}  on $[T_{n-3}, T_{n-1}]$ with terminal condition $V(T_{n-1},\cdot,\cdot)$. We emphasize that  $\pa_xV(T_{n-1},\cdot,\cdot)$ has the above uniform regularity with the same constants $L^V_{xx}, C^\mu_2$, then we may apply Proposition \ref{prop:pamuV}  with the same $\d$ and obtain a classical solution $V$ on $[T_{n-3}, T_{n-1}]$ with the additional regularities specified in Proposition \ref{prop:pamuV}-(iii)(iv). Clearly this extends the classical solution of the master equation to $[T_{n-3}, T_n]$. We emphasize again that, while the bound of $\pa_{x\mu} V(t,\cd)$ may become larger for $t\in [T_{n-3}, T_{n-2}]$ because the $C^\mu_1$ in \reff{W1Lip} now depends on $\|\pa_{x\mu} V(T_{n-1}, \cd)\|_{L^\infty}$ instead of $\|\pa_{x\mu} V(T_{n}, \cd)\|_{L^\infty}$, by the global a priori estimates in Theorems \ref{thm-antidisplacement} and \ref{thm-Lipschitz2} we see that $\pa_xV(t,\cd)$ corresponds to the same $L_{xx}^V$ and $C^\mu_2$ for all $t\in [T_{n-3}, T_{n}]$. This enables us to consider the master equation \reff{master}  on $[T_{n-4}, T_{n-2}]$ with terminal condition $V(T_{n-2},\cdot,\cdot)$,  and then we obtain a classical solution on $[T_{n-4}, T_n]$ with the desired uniform estimates and additional regularities. 
 
Repeat the arguments backwardly in time, we may construct a classical solution $V$ for the original master equation \reff{master} on $[0, T]$ with terminal condition $G$. Moreover, since this procedure is repeated only $n$ times, by applying \reff{W1Lip} repeatedly we see that \reff{W1Lip} indeed holds true on $[0, T]$.
\qed

We conclude the paper by providing an example which satisfies all the assumptions in Theorem \ref{thm:wellposedness}.
We emphasize that there is no smallness assumption imposed here.

\begin{eg}\label{eg:allassum}
For simplicity let $d=1$. Fix positive constants $0<\underline \a\leq \overline\a$ and $0<\underline \gamma<\overline\gamma$   with $\overline\gamma>1$, and fix $(\l_1,\l_2,\l_3)$ satisfying the requirements in \reff{D4}. Set $\ol L^A:=1$ and let $M_0$ be a large number which will be specified later. Assume 

(i) $G$ satisfies Assumption \ref{assum-regG} with 
\begin{equation}\label{eq:paxxGkappa}
-\overline\a M_0\leq \pa_{xx}G(x,\mu)\leq -\underline\a M_0\quad {\rm{on}}\q \mathbb R\times\mathcal{P}_2(\dbR);
\end{equation}

(ii)  $H$ satisfies Assumption \ref{assum-regH1} with $A_0:= M_0^3 > L^{H_0}_2$ in  \eqref{eq:Hform}, and
\begin{equation}\label{paxxH0exp}
  \underline \g[A_0-L_2^{H_0}]\leq \pa_{xx}H_0(x,\mu,p)\leq \overline\g[A_0-L_2^{H_0}]\q {\rm{on}}\q \mathbb R\times\mathcal{P}_2(\dbR)\times\dbR.
\end{equation}
Then, for $M_0$ large enough, which may depend on $\ul \a, \ol \a, \ul \g, \ol \g$, $(\l_1, \l_2, \l_3)$, and $L_2^G$,  $L_2^{H_0}$, one can choose appropriate $\l_0$ such that all the conditions in Theorem \ref{thm:wellposedness} hold true.
\end{eg}
\proof We first emphasize that \reff{eq:paxxGkappa} and \reff{paxxH0exp} involve only $\pa_{xx} G$ and $\pa_{xx} H_0$. Note that the parameters $L_2^G$,  $L_2^{H_0}$, which $M_0$ will depend on, do not involve these derivatives. So it is rather easy to construct $G$ and $H_0$ satisfy both  Assumptions \ref{assum-regG}, \ref{assum-regH1}, and  \reff{eq:paxxGkappa}, \reff{paxxH0exp} with arbitrarily large $M_0$. Moreover, recall \reff{Gbound} and \reff{H0bound1}, by \reff{eq:paxxGkappa} and \reff{paxxH0exp} it is clear that
\bea
\label{newLGxx}
\overline L_{xx}^G=\overline\a M_0, \q \underline L_{xx}^{H_0}=\underline \g[A_0-L_2^{H_0}],\q \overline L_{xx}^{H_0}=\overline\g[A_0-L_2^{H_0}].
\eea
Then the $\th_3$ in \reff{LA} and $L^u_{xx}(\th_3)$ in \reff{pauest1} become: recalling $\ol L^A=1$, 
\bea
\label{newtheta3}
\left.\ba{c}
\dis \theta_3 := 1+L_2^{H_0}\big[1+\overline\a M_0+\sqrt{(1+\overline\a M_0)^2-1}\big],\ms\\
\dis L_{xx}^V = L_{xx}^u(\theta_3) := \frac{2\overline\a M_0(\theta_3-1)}{\theta_3-1-L_2^{H_0}+\sqrt{(\theta_3-1-L_2^{H_0})^2-2L_2^{H_0}\overline \a M_0(\theta_3-1)}}.
\ea\right.
\eea
We now show that the following $\l_0$ satisfies all the requirements:
\bea
\label{lambda0000}
\l_0=\frac{\bar\gamma^2[1+L_{xx}^V]^2-8\l_3}{4\underline \gamma}+1.
\eea

First, by the choice of $\l_0$, it is obvious that 
$
\l_0>\frac{\bar\gamma^2[1+L_{xx}^V]^2-8\l_3}{4\underline \gamma},
$
which verifies \eqref{theta}. 

Next, let $O(M)$ denote a generic positive function of $M$ such that ${O(M)\over M}$ is bounded both from above and away from $0$. Then we see that 
 \begin{equation}\label{newtheta3LxxV}
\theta_3=O(M_0), \q L_{xx}^V=O(M_0),\q \l_0 = O(M_0^2).
\end{equation}
By \reff{anti} we have
\beaa
(AntiMon)^{\vec \l}_\xi U(\eta, \eta) \le \Big[- \l_0 \ul \a M_0 + \ol\a^2 M_0^2 -\l_3 \Big] \dbE[|\eta|^2] + \Big[|\l_1| L^G_2 + \l_2|L^G_2|^2\Big] \big|\dbE[\eta]\big|^2.
\eeaa
Since $\l_0 M_0 = O(M_0^3)$, it is clear that $G$ is $\vec{\l}$-anti-monotone when $M_0$ is large enough.

Moreover,  since $d=1$, we have $\underline \kappa(A_0)=\overline \kappa(A_0)=\underline \kappa'(A_0)=A_0$ and $L^{A_0}=1\leq \overline L^A=1$. Recall \reff{theta} and \reff{A}. When $M_0$ is large, it is clear that $1-\th_1$ is uniformly away from $0$ and then it follows from \reff{newtheta3LxxV} that $\underline \k(A_1^{-1}A_2) = O(M_0^2)$.   Thus, since $A_0=M_0^3$ and  $\overline \gamma>1$, for $M_0$ sufficiently large we have the following inequalities which verify \eqref{kA0}:
\[
A_0\geq 1+\underline \k(A_1^{-1}A_2)L_2^{H_0},\quad A_0\geq\theta_3\quad\text{and}\quad A_0+L_2^{H_0}\leq\overline\g[A_0-L_2^{H_0}].
\]

Finally, since  $\overline L_{xx}^G=\overline\a M_0$, it is clear that $2\overline L_{xx}^G[A_0-1]\geq\overline\g[A_0-L_2^{H_0}]$ for $M_0$ large enough. Then \reff{paxxH0exp} implies  \eqref{LH0}.
\qed

\end{document}